\documentclass[aos, preprint]{imsart}

\RequirePackage[OT1]{fontenc}
\RequirePackage{amsthm,amsmath}
\numberwithin{equation}{section}
\theoremstyle{plain}
\newtheorem{thm}{Theorem}

\newtheorem{lem}{Lemma}

\newcommand{\cX}{{\cal X}}

\newcommand{\cZ}{{\cal Z}}

\newcommand{\wht}{\widehat}
\newcommand{\wtd}{\widetilde}

\newcommand{\wbox}{\sqcap\llap{$\sqcup$}}

\begin{document}
\begin{frontmatter}

% "Title of the paper"
\title{THE MULTI-ARMED BANDIT PROBLEM: AN EFFICIENT NON-PARAMETRIC SOLUTION}

% indicate corresponding author with \corref{}
% \author{\fnms{John} \snm{Smith}\corref{}\ead[label=e1]{smith@foo.com}\thanksref{t1}}
% \thankstext{t1}{Thanks to somebody} 
% \address{line 1\\ line 2\\ printead{e1}}
% \affiliation{Some University}

\author{Hock Peng Chan\ead[label=e1]{stachp@nus.edu.sg}\thanksref{t1}}
\thankstext{t1}{Supported by MOE grant number R-155-000-158-112}
\affiliation{National University of Singapore}

\address{Department of Statistics and Applied Probability \\
Block S16, Level 7, 6 Science Drive 2 \\
Faculty of Science \\
National University of Singapore \\
Singapore 117546}

%\and
%\author{\fnms{Hao} \snm{Chen}\ead[label=e2]{gwalther@stat.stanford.edu}\thanksref{t2}}
%\thankstext{t2}{?}
%\address{Statistics Department}
%\affiliation{University of California at Davis}
%\and
%\author{\fnms{???} \snm{???}\ead[label=e2]{???}}
%\address{\printead{e2}}
%\affiliation{???}

\begin{abstract}
Lai and Robbins (1985) and Lai (1987) provided efficient parametric solutions to the multi-armed bandit problem,
showing that arm allocation via upper confidence bounds (UCB) achieves minimum regret.
These bounds are constructed from the Kullback-Leibler information of the reward distributions,
estimated from specified parametric families.
In recent years there has been renewed interest in the multi-armed bandit problem due to new applications in machine learning algorithms and data analytics. 
Non-parametric arm allocation procedures like $\epsilon$-greedy, Boltzmann exploration and BESA were studied,
and modified versions of the UCB procedure were also analyzed under non-parametric settings.
However unlike UCB these non-parametric procedures are not efficient under general parametric settings.
In this paper we propose efficient non-parametric procedures. 
\end{abstract}


%\kwd[; secondary ]{60J22}
%\kwd{60K35}
%\end{keyword}

\end{frontmatter}

\section{Introduction}

Lai and Robbins (1985) provided an asymptotic lower bound for the regret in the multi-armed bandit problem, 
and proposed an index strategy that is efficient,
that is it achieves this bound. 
Lai (1987) showed that allocation to the arm having the highest upper confidence bound (UCB), 
constructed from the Kullback-Leibler (KL) information between the estimated reward distributions of the arms, 
is efficient when the distributions belong to a specified exponential family.
Agrawal (1995) proposed a modified UCB procedure that is efficient despite not having to know in advance the total sample size.
Capp\'{e}, Garivier, Maillard, Munos and Stoltz (2013) provided explicit,
non-asymptotic bounds on the regret of a KL-UCB procedure that is efficient on a larger class of distribution families.

Burnetas and Kalehakis (1996) extended UCB to multi-parameter families,
almost showing efficiency in the natural setting of normal rewards with unequal variances.
Yakowitz and Lowe (1991) proposed non-parametric procedures that do not make use of KL-information,
suggesting logarithmic and polynomial rates of regret under finite exponential moment and moment conditions respectively.

Auer, Cesa-Bianchi and Fischer (2002) proposed a UCB1 procedure that achieves logarithmic regret when the reward distributions are supported on [0,1].
They also studied the $\epsilon$-greedy algorithm of Sutton and Barto (1998) and provided finite-time upper bounds of its regret.
Both UCB1 and $\epsilon$-greedy are non-parametric in their applications and,
unlike UCB-Lai or UCB-Agrawal,
are not expected to be efficient under a general exponential family setting.
Other non-parametric methods that have been proposed include reinforcement comparison, 
Boltzmann exploration (Sutton and Barto, 1998) and pursuit (Thathacher and Sastry, 1985).
Kuleshov and Precup (2014) provided numerical comparisons between UCB and these methods.
For a description of applications to recommender systems and clinical trials,
see Shivaswamy and Joachims (2012).
Burtini, Loeppky and Lawrence (2015) provided a comprehensive survey of the methods, results and applications of the multi-armed bandit problem,
developed over the past thirty years. 

A strong competitor to UCB under the parametric setting is the Bayesian method,
see for example Fabius and van Zwet (1970) and Berry (1972). 
There is also a well-developed literature on optimization under an infinite-time discounted window setting,
in which allocation is to the arm maximizing a dynamic allocation (or Gittins) index,
see the seminal papers Gittins (1979) and Gittins and Jones (1979), 
and also Berry and Fristedt (1985), Chang and Lai (1987), Brezzi and Lai (2002).
Recently there has been renewed interest in the Bayesian method due to the developments of UCB-Bayes [see Kaufmann, Capp\'{e} and Garivier (2012)] 
and Thompson sampling [see for example Korda, Kaufmann and Munos (2013)].

In this paper we propose an arm allocation procedure subsample-mean comparison (SSMC),
that though non-parametric,
is nevertheless efficient when the reward distributions are from an {\it unspecified} one-dimensional exponential family.
It achieves this by comparing subsample means of the leading arm with the sample means of its competitors.
It is empirical in its approach, 
using more informative subsample means rather than full-sample means alone,
for better decision-making. 
The subsampling strategy was first employed by Baransi, Maillard and Mannor (2014) in their best empirical sampled average (BESA) procedure.
However there are key differences in their implementation of subsampling from ours,
as will be elaborated in Section~2.2.
Though efficiency has been attained for various one-dimensional exponential families by say UCB-Agrawal or KL-UCB,
SSMC is the first to achieve efficiency without having to know the specific distribution family. 
In addition we propose in Section 2.4 a related subsample-$t$ comparison (SSTC) procedure,
applying $t$-statistic comparisons in place of mean comparisons,
that is efficient for normal distributions with unknown and unequal variances.

The layout of the paper is as follows.
In Section 2 we describe the subsample comparison strategy for allocating arms.
In Section~3 we show that the strategy is efficient for exponential families,
including the setting of normal rewards with unknown and unequal variances.
In Section~4 we show logarthmic regret for Markovian rewards.
In Section~5 we provide numerical comparisons against existing methods.
In Section 6 we provide a concluding discussion.
In Section~7 we prove the results of Sections~3 and 4.

\section{Subsample comparisons}

Let $Y_{k1}, Y_{k2}, \ldots$,
$1 \leq k \leq K$,
be the observations (or rewards) from a population (or arm) $\Pi_k$.
We assume here and in Section 3 that the rewards are independent and identically distributed (i.i.d.) within each arm.
We extend to Markovian rewards in Section~4.
Let $\mu_k = E Y_{kt}$ and $\mu_* = \max_{1 \leq k \leq K} \mu_k$.

Consider a sequential procedure for selecting the population to be sampled,
with the decision based on past rewards.
Let $N_k$ be the number of observations from $\Pi_k$ when there are $N$ total observations,
hence $N = \sum_{k=1}^K N_k$.
The objective is to minimize the {\it regret}
$$R_N := \sum_{k=1}^K (\mu_*-\mu_k) EN_k.
$$

The Kullback-Leibler information number between two densities $f$ and $g$,
with respect to a common ($\sigma$-finite) measure,
is 
\begin{equation} \label{Dfg}
D(f|g) = E_f [ \log \tfrac{f(Y)}{g(Y)}],
\end{equation}
where $E_f$ denotes expectation with respect to $Y \sim f$. 
An arm allocation procedure is said to be uniformly good if
\begin{equation} \label{fast}
R_N = o(N^{\epsilon}) \mbox{ for all } \epsilon > 0,
\end{equation}
over all reward distributions lying within a specified parametric family. 

Let $f_k$ be the density of $Y_{kt}$ and let $f_* = f_k$ for $k$ such that $\mu_k = \mu_*$ (assuming $f_*$ is unique).
The celebrated result of Lai and Robbins (1985) is that under (\ref{fast}) and additional regularity conditions,
\begin{equation} \label{RN}
\liminf_{N \rightarrow \infty} \frac{R_N}{\log N} \geq \sum_{k: \mu_k < \mu_*} \frac{\mu_*-\mu_k}{D(f_k|f_*)}.
\end{equation}
Lai and Robbins (1985) and Lai (1987) went on to propose arm allocation procedures that have regrets achieving the lower bound in (\ref{RN}),
and are hence {\it efficient}.

\subsection{Review of existing methods}

In the setting of normal rewards with unit variances,
UCB-Lai can be described as the selection,
for sampling,
$\Pi_k$ maximizing
\begin{equation} \label{UCBLai}
\bar Y_{k n_k} + \sqrt{\tfrac{2 \log(N/n)}{n}},
\end{equation}
where $\bar Y_{kt} = \frac{1}{t} \sum_{u=1}^t Y_{ku}$,
$n$ is the current number of observations from the $K$ populations,
and $n_k$ is the current number of observations from $\Pi_k$.
Agrawal (1995) proposed a modified version of UCB-Lai that does not involve the total sample size $N$, 
with the selection instead of the population $\Pi_k$ maximizing
\begin{equation} \label{Agr}
\bar Y_{k n_k} + \sqrt{\tfrac{2 (\log n+\log \log n+b_n)}{n_k}}, 
\end{equation}
with $b_n \rightarrow \infty$ and $b_n = o(\log n)$.
Efficiency holds for (\ref{UCBLai}) and (\ref{Agr}),
and there are corresponding versions of (\ref{UCBLai}) and (\ref{Agr}) that are efficient for other one-parameter exponential families.
Capp\'{e} et al. (2013) proposed a more general KL-UCB procedure that is also efficient for distributions with given finite support.

Auer, Cesa-Bianchi and Fischer (2002) simplified UCB-Agrawal to UCB1, 
proposing that $\Pi_k$ maximizing
\begin{equation} \label{ucb1}
\bar Y_{k n_k} + \sqrt{\tfrac{2 \log n}{n_k}}
\end{equation}
be selected. 
They showed that under UCB1,
logarithmic regret $R_N = O(\log N)$ is achieved when the reward distributions are supported on [0,1].
In the setting of normal rewards with unequal and unknown variances,
Auer et al. suggested applying a variant of UCB1 which they called UCB1-Normal,
and showed logarithmic regret.
Under UCB1-Normal,
an observation is taken from any population $\Pi_k$ with $n_k < 8 \log n$.
If such a population does not exist,
then an observation is taken from $\Pi_k$ maximizing
$$\bar Y_{k n_k} + 4 \wht \sigma_{kn_k} \sqrt{\tfrac{\log n}{n_k}},
$$
where $\wht \sigma_{kt}^2 = \frac{1}{t-1} \sum_{u=1}^{t} (Y_{ku}-\bar Y_{kt})^2$.

Auer et al. provided an excellent study of various non-parametric arm allocation procedures,
for example the $\epsilon$-greedy procedure proposed by Sutton and Barto (1998),
in which an observation is taken from the population with the largest sample mean with probability $1-\epsilon$,
and randomly with probability $\epsilon$.
Auer et al. suggested replacing the fixed $\epsilon$ at every stage by a stage-dependent 
\begin{equation} \label{epsn}
\epsilon_n = \min(1,\tfrac{cK}{d^2n}),
\end{equation}
with $c$ user-specified and $0 < d \leq \min_{k: \mu_k < \mu^*} (\mu_*-\mu_k)$.
They showed that if $c>5$,
then logarithmic regret is achieved for reward distributions supported on $[0,1]$.
A more recent numerical study by Kuleshov and Precup (2014) considered additional non-parametric procedures,
for example Boltzmann exploration in which an observation is taken from $\Pi_k$ with probability proportional to $e^{\bar Y_{kn_k}/\tau}$, 
for some $\tau > 0$.

\subsection{Subsample-mean comparisons}

A common characteristic of the procedures described in Section 2.1 is that allocation is based solely on a comparison of the sample means $\bar Y_{kn_k}$,
with the exception of UCB1-Normal in which $\wht \sigma_{kn_k}$ is also utilized.
As we shall illustrate in Section 2.3,
we can utilize subsample-mean information from the leading arm to estimate the confidence bounds for selecting from the other arms.
In contrast UCB-based procedures like KL-UCB discard subsample information and rely on parametric information to estimate these bounds.
Even though subsample-mean and KL-UCB are both efficient for exponential families,
the advantage of subsample-mean is that the underlying family need not be specified.

In SSMC a leader is chosen in each round of play to compete against all the other arms.
Let $r$ denote the round number.
In round 1,
we sample all $K$ arms.
In round $r$ for $r>1$,
we set up a challenge between the leading arm (to be defined below) and each of the other arms.
An arm is sampled only if it wins all its challenges in that round.
Hence for round $r>1$ we sample either the leading arm or a non-empty subset of the challengers.  
Let $n(=n^r)$ be the total number of observations from all $K$ arms at the beginning of round~$r$,
let $n_k(=n_k^r)$ be the corresponding number from $\Pi_k$.
Hence $n_k^1=0$ and $n_k^2=1$ for all $k$,
and $K+(r-2) \leq n^r \leq K+(K-1)(r-2)$ for $r \geq 2$.

Let $c_n$ be a non-negative monotone increasing sampling threshold in SSMC and SSTC,
with
\begin{equation} \label{cnn}
c_n = o(\log n) \mbox{ and } \tfrac{c_n}{\log \log n} \rightarrow \infty \mbox{ as } n \rightarrow \infty.
\end{equation}
For example in our implementation of SSMC and SSTC in Section 5,
we select $c_n=(\log n)^{\frac{1}{2}}$. 
An explanation of why (\ref{cnn}) is required for efficiency of SSMC is given in the beginning of Section 7.1.
Let $\bar Y_{k,t:u} = \frac{1}{u-t+1} \sum_{v=t}^u Y_{kv}$,
hence $\bar Y_{kt} = \bar Y_{k,1:t}$.

\smallskip
\underline{Subsample-mean comparison (SSMC)}
\begin{enumerate}
\item $r=1$. 
Sample each $\Pi_k$ exactly once.

\item $r=2,3,\ldots$.
\begin{enumerate}
\item Let the leader $\zeta(=\zeta^r)$ be the population with the most observations,
with ties resolved by (in order):

\begin{enumerate}
\item the population with the larger sample mean,
\item the leader of the previous round,
\item randomization.
\end{enumerate}

\item For all $k \neq \zeta$ set up a challenge between $\Pi_{\zeta}$ and $\Pi_k$ in the following manner.

\begin{enumerate}
\item If $n_k = n_{\zeta}$, 
then $\Pi_k$ loses the challenge automatically.

\item If $n_k < n_{\zeta}$ and $n_k < c_n$, 
then $\Pi_k$ wins the challenge automatically.

\item If $c_n \leq n_k < n_{\zeta}$,
then $\Pi_k$ wins the challenge when
\begin{equation} \label{meancompare}
\bar Y_{kn_k} \geq \bar Y_{\zeta,t:(t+n_k-1)} \mbox{ for some } 1 \leq t \leq n_{\zeta}-n_k+1.
\end{equation}
\end{enumerate}

\item For all $k \neq \zeta$, 
sample from $\Pi_k$ if $\Pi_k$ wins its challenge against $\Pi_{\zeta}$. 
Sample from $\Pi_{\zeta}$ if $\Pi_{\zeta}$ wins all its challenges.
Hence either $\Pi_{\zeta}$ is sampled,
or a non-empty subset of $\{ \Pi_k: k \neq \zeta \}$ is sampled.
\end{enumerate}
\end{enumerate}

\smallskip
SSMC may recommend more than one populations to be sampled in a single round when $K>2$.
In the event that $n^r < N < n^{r+1}$ for some $r$,
we select $N-n^r$ populations randomly from among the $n^{r+1}-n^r$ recommended by SSMC in the $r$th round,
in order to make up exactly $N$ observations.

If $\Pi_{\zeta}$ wins all its challenges,
then $\zeta$ and $(n_k: k \neq \zeta)$ are unchanged, 
and in the next round it suffices to perform the comparison in (\ref{meancompare}) at the largest $t$ instead of at every $t$.
The computational cost is thus $O(1)$.
The computational cost is $O(r)$ if at least one $k \neq \zeta$ wins its challenge.
Hence when there is only one optimal arm and SSMC achieves logarithmic regret,
the total computational cost is $O(r \log r)$ for $r$ rounds of the algorithm.

In step 2(b)ii. we force the exploration of arms with less than $c_n$ rewards.
By (\ref{cnn}) we select $c_n$ small compared to $\log n$,
so that the cost of such forced explorations is asymptotically negligible.
In contrast the forced exploration in the greedy algorithm (\ref{epsn}) is more substantial,
of order $\log n$ for $n$ rewards.

BESA,
proposed by Baransi, Maillard and Mannor (2014),
also applies subsample-mean comparisons.
We describe BESA for $K=2$ below,
noting that tournament-style elimination is applied for $K>2$.
Unlike SSMC,
exactly one population is sampled in each round $r>1$ even when $K>2$.

\smallskip
\underline{Best Empirical Sampled Average (BESA)}

\begin{enumerate}
\item $r=1$. Sample both $\Pi_1$ and $\Pi_2$. 

\item $r=2,3,\ldots$.
\begin{enumerate}
\item Let the leader $\zeta$ be the population with more observations,
and let $k \neq \zeta$.

\item Sample randomly without replacement $n_k$ of the $n_{\zeta}$ observations from $\Pi_{\zeta}$, 
and let $\bar Y_{\zeta n_k}^*$ be the mean of the $n_k$ observations.

\item If $\bar Y_{k n_k} \geq \bar Y_{\zeta n_k}^*$,
then sample from $\Pi_k$.
Otherwise sample from $\Pi_{\zeta}$.
\end{enumerate}
\end{enumerate}

\smallskip
As can be seen from the descriptions of SSMC and BESA,
the mechanism of choosing the arm to be played in SSMC clearly promotes exploration of non-leading arms,
relative to BESA.
Whereas Baransi et al. demonstrated logarithmic regret of BESA for rewards bounded on [0,1] (though BESA can of course be applied on more general settings but with no such guarantees),
we show in Section 3 that SSMC is able to extend BESA's subsampling idea to achieve asymptotic optimality,
that is efficiency,
on a wider set of distributions. 
Tables 4 and 5 in Section 5 show that SSMC controls the oversampling of inferior arms better relative to BESA,
due to its added explorations.

\subsection{Comparison of SSMC with UCB methods}

Lai and Robbins (1985) proposed a UCB strategy in which the arms take turns to challenge a leader with order $n$ observations. 
Let us restrict to the setting of exponential families.
Denote the leader by $\zeta$ and the challenger by $k$.
Lai and Robbins proposed,
in their (3.1),
upper confidence bounds $U_{kt}^n = U_k^n(Y_{k1}, \ldots, Y_{kt})$ satisfying 
$$P(\min_{1 \leq t \leq n} U_{kt}^n \geq \mu_k-\epsilon) = 1-o(n^{-1}) \mbox{ for all } \epsilon>0.
$$
The decision is to sample from arm $k$ if
$$U_{kn_k}^n \geq \bar Y_{\zeta n_{\zeta}} (\doteq \mu_{\zeta}),
$$
otherwise arm $\zeta$ is sampled.
By doing this we ensure that if $\mu_k > \mu_{\zeta}$,
then the probability that arm~$k$ is sampled is $1-o(n^{-1})$.

We next consider SSMC. 
Let $L_{\zeta n_k} = \min_{1 \leq t \leq n_{\zeta}-n_k+1} \bar Y_{\zeta,t:(t+n_k-1)}$.
Since $n_{\zeta}$ is of order $n$,
it follows that if $\mu_k > \mu_{\zeta}$,
then as $Y_{kt}$ is stochastically larger than $Y_{\zeta t}$, 
$$P(L_{\zeta n_k} \leq \bar Y_{k n_k}) = 1-o(n^{-1}).
$$
In SSMC we sample from arm $k$ if $L_{\zeta n_k} \leq \bar Y_{k n_k}$,
ensuring, 
as in Lai and Robbins, 
that an optimal arm is sampled with probability $1-o(n^{-1})$ when the leading arm is inferior.

In summary SSMC differs from UCB in that it compares $\bar Y_{kn_k}$ against a lower confidence bound $L_{\zeta n_k}$ of the leading arm,
computed from subsample-means instead of parametrically.
Nevertheless the critical values that SSMC and UCB-based methods employ for allocating arms are asymptotically the same, 
as we shall next show.

For simplicity let us consider unit variance normal densities with $K=2$.
Consider firstly unbalanced sample sizes with say $n_2=O(\log n)$ and note,
see Appendix A,
that  
\begin{equation} \label{min1}
\min_{1 \leq t \leq n_1-n_2+1} \bar Y_{1,t:(t+n_2-1)} = \mu_1 - [1+o_p(1)] \sqrt{\tfrac{2 \log n}{n_2}}.
\end{equation}
Hence arm 2 winning the challenge requires
\begin{equation} \label{compare}
\bar Y_{2n_2} \geq \mu_1 - [1+o_p(1)] \sqrt{\tfrac{2 \log n}{n_2}}.
\end{equation}
By (\ref{Agr}) and (\ref{ucb1}),
UCB-Agrawal, KL-UCB and UCB1 also select arm 2 when (\ref{compare}) holds,
since $\bar Y_{1n_1}+\sqrt{\frac{2 \log n}{n_1}} = \mu_1+o_p(1)$.
Hence what SSMC does is to estimate the critical value $\mu_1 - [1+o_p(1)] \sqrt{\frac{2 \log n}{n_2}}$, 
empirically by using the minimum of the running averages $\bar Y_{1,t:(t+n_2-1)}$.
In the case of $n_1,n_2$ both large compared to $\log n$,
$\sqrt{\frac{2 \log n}{n_1}}+\sqrt{\frac{2 \log n}{n_2}} \rightarrow 0$,
and SSMC, UCB-Agrawal, KL-UCB and UCB1 essentially select the population with the larger sample mean.

\subsection{Subsample-$t$ comparisons}

For efficiency outside one-parameter exponential families,
we need to work with test statistics beyond sample means.
For example to achieve efficiency for normal rewards with unknown and unequal variances,
the analogue of mean comparisons is $t$-statistic comparisons
$$\frac{\bar Y_{kn_k}-\mu_{\zeta}}{\hat \sigma_{kn_k}} \geq \frac{\bar Y_{\zeta,t:(t+n_k-1)}-\mu_{\zeta}}{\hat \sigma_{\zeta,t:(t+n_k-1)}},
$$
where $\wht \sigma^2_{k,t:u} = \frac{1}{u-t} \sum_{v=t}^u (Y_{kv}-\bar Y_{k,t:u})^2$ and $\wht \sigma_{kt} = \wht \sigma_{k,1:t}$.
Since $\mu_{\zeta}$ is unknown,
we estimate it by $\bar Y_{\zeta n_{\zeta}}$.

\smallskip
\underline{Subsample-$t$ comparison (SSTC)}

Proceed as in SSMC,
with step 2(b)iii.$'$ below replacing step 2(b)iii.

iii.$'$ If $c_n \leq n_k < n_{\zeta}$,
then $\Pi_k$ wins the challenge when either $\bar Y_{k n_k} \geq \bar Y_{\zeta n_{\zeta}}$ or
\begin{equation} \label{comparet}
\frac{\bar Y_{kn_k}-\bar Y_{\zeta n_{\zeta}}}{\hat \sigma_{k n_k}} \geq \frac{\bar Y_{\zeta,t:(t+n_k-1)}-\bar Y_{\zeta n_{\zeta}}}{\hat \sigma_{\zeta,t:(t+n_k-1)}}
\mbox{ for some } 1 \leq t \leq n_{\zeta}-n_k+1.
\end{equation}

\smallskip
As in SSMC only $O(r \log r)$ computations are needed for $r$ rounds when there is only one optimal arm and the regret is logarithmic.
This is because it suffices to record the range of $\bar Y_{\zeta n_{\zeta}}$ that satisfies (\ref{comparet}) for each $k \neq \zeta$,
and the actual value of $\bar Y_{\zeta n_{\zeta}}$.
The updating of these requires $O(1)$ computations when both $\zeta$ and $(n_k: k \neq \zeta)$ are unchanged. 

\section{Efficiency}

Consider firstly an exponential family of density functions
\begin{equation} \label{expo}
f(x;\theta) = e^{\theta x-\psi(\theta)} f(x;0), 
\ \theta \in \Theta,
\end{equation}
with respect to some measure $\nu$,
where $\psi(\theta) = \log [\int e^{\theta x} f(x;0) \nu(dx)]$ is the log moment generating function and $\Theta = \{ \theta: \psi(\theta) < \infty \}$.
For example the Bernoulli family satisfies (\ref{expo}) with $\nu$ the counting measure on $\{ 0,1 \}$ and $f(0;0)= f(1;0)=\frac{1}{2}$.
The family of normal densities with variance $\sigma^2$ satisfies (\ref{expo}) with $\nu$ the Lebesgue measure and $f(x;0) = \frac{1}{\sigma \sqrt{2 \pi}} e^{-x^2/(2 \sigma^2)}$.

Let $f_k = f(\cdot;\theta_k)$ for some $\theta_k \in \Theta$,
$1 \leq k \leq K$.
Let $\theta_* = \max_{1 \leq k \leq K} \theta_k$ and $f_* = f(\cdot;\theta_*)$.
By (\ref{Dfg}) and (\ref{expo}),
the KL-information in (\ref{RN}),
\begin{eqnarray*}
D(f_k|f_*) & = & \int \{ (\theta_k-\theta_*)x-[\psi(\theta_k)-\psi(\theta_*)] \} f(x; \theta_k) \nu(dx) \cr
& = & (\theta_k-\theta_*) \mu_k - [\psi(\theta_k)-\psi(\theta_*)] = I_*(\mu_k),
\end{eqnarray*}
where $I_*$ is the large deviations rate function of $f_*$.
Let $\Xi = \{ \ell: \mu_{\ell} = \mu_* \}$ be the set of optimal arms.

\begin{thm} \label{thm1}
For the exponential family {\rm (\ref{expo})},
SSMC satisfies
\begin{equation} \label{1eq}
\limsup_{r \rightarrow \infty} \frac{En_k^r}{\log r} \leq \frac{1}{D(f_k|f_*)}, 
\quad k \not\in \Xi,
\end{equation}
and is thus efficient.
\end{thm}

UCB-Agrawal and KL-UCB are efficient as well for (\ref{expo}),
see Agrawal (1995) and Capp\'{e} et al. (2013),
SSMC is unique in that it achieves efficiency by being adaptive to the exponential family,
whereas UCB-Agrawal and KL-UCB achieve efficiency by having selection procedures that are specific to the exponential family.
On the other hand UCB-based methods require less storage space, 
and more informative finite-time bounds have been obtained. 
Specifically for UCB-based methods in exponential families we need only store the sample mean for each arm,
and the numerical complexity is of the same order as the sample size.
For SSMC as given in Section 2.3,
all observations are stored (more of this in Section 6) and the numerical complexity for a sample of size $N$ is $N \log N$ when we have efficiency and exactly one optimal arm.

We next consider normal rewards with unequal and unknown variances, 
that is with densities 
\begin{equation} \label{fnormal}
f(x;\mu,\sigma^2) = \tfrac{1}{\sigma \sqrt{2 \pi}} e^{-\frac{(x-\mu)^2}{2 \sigma^2}},
\end{equation}
with respect to Lebesgue measure.
Let $M(g) = \frac{1}{2} \log(1+g^2)$.
Burnetas and Katehakis (1996) showed that if $f_k = f(\cdot;\mu_k,\sigma_k^2)$, 
then under uniformly fast convergence and additional regularity conditions,
an arm allocation procedure must have regret $R_N$ satisfying
$$\liminf_{N \rightarrow \infty} \frac{R_N}{\log N} \geq \sum_{k: \mu_k < \mu_*} \frac{\mu_*-\mu_k}{M(\frac{\mu_*-\mu_k}{\sigma_k})}. 
$$
They proposed an extension of UCB-Lai but needed the verification of a technical condition to show efficiency.
In the case of UCB1-Normal, 
logarithmic regret also depended on tail bounds of the $\chi^2$- and $t$-distributions that were only shown to hold numerically by Auer et al. (2002). 
In Theorem~\ref{thm2} we show that SSTC achieves efficiency. 

\begin{thm} \label{thm2}
For normal densities {\rm (\ref{fnormal})} with unequal and unknown variances,
SSTC satisfies
$$\limsup_{r \rightarrow \infty} \frac{E n_k^r}{\log r} \leq \frac{1}{M(\frac{\mu_*-\mu_k}{\sigma_k})},
\quad k \not\in \Xi,
$$
and is thus efficient.
\end{thm}

\section{Logarithmic regret}

We show here that logarithmic regret can be achieved by SSMC under Markovian assumptions. 
This is possible because in SSMC we compare blocks of observations that retain the Markovian structure.

For $1 \leq k \leq K$,
let $X_{k1}, X_{k2}, \ldots$ be a potentially unobserved ${\cal X}$-valued Markov chain,
with $\sigma$-field ${\cal A}$ and transition kernel
\begin{equation} \label{PYB}
P_k(x,A) = P (X_{kt} \in A| X_{k,t-1} = x), \ x \in {\cal X}, A \in {\cal A}.
\end{equation}
We shall assume for convenience that $(X_{kt})_{t \geq 1}$ is stationary.
Let $Y_{k1}, Y_{k2}, \ldots$ be real-valued and conditionally independent given $(X_{kt})_{t \geq 1}$,
and having conditional densities $\{ f_k(\cdot|x): 1 \leq k \leq K, x \in {\cal X} \}$,
with respect to some measure $\nu$, 
such that
$$P ( Y_{kt} \in B| X_{k1}=x_1, X_{k2} = x_2, \cdots ) = \int_B f_k(y|x_t) \nu(dy).
$$
We assume that the $K$ Markov chains are independent, 
and that the following Doeblin-type condition holds.

\smallskip
\noindent (C1) For $1 \leq k \leq K$,
there exists a non-trival measure $\lambda_k$ on $({\cal X},{\cal A})$ such that
$$P_k(x,A) \geq \lambda_k(A), \
x \in {\cal X}, A \in {\cal A}.
$$

\smallskip
As before let $\mu_k = E Y_{kt}$,
$\mu_* = \max_{1 \leq k \leq K} \mu_k$ and the regret 
$$R_N = \sum_{k: \mu_k < \mu_*} (\mu_*-\mu_k) EN_k.
$$
In addition to (C1) we assume the following sample mean large deviations.  

\smallskip
\noindent (C2) For any $\epsilon > 0$,
there exists $b(=b_{\epsilon}) > 0$ and $Q(=Q_{\epsilon}) > 0$ such that for $1 \leq k \leq K$ and $t \geq 1$,
\begin{equation} \label{YQ}
P ( |\bar Y_{kt} - \mu_k| \geq \epsilon ) \leq Qe^{-t b}. 
\end{equation}

\smallskip
\noindent (C3) For $k$ such that $\mu_k < \mu_*$ and $\ell$ such that $\mu_{\ell} = \mu_*$,
there exists $b_1 > 0$, 
$Q_1>0$ and $t_1 \geq 1$ such that for $\omega \leq \mu_k$ and $t \geq t_1$,
\begin{equation} \label{Q1b}
P(\bar Y_{\ell t} < \omega) \leq Q_1 e^{-t b_1} P(\bar Y_{kt} < \omega).
\end{equation}

\begin{thm} \label{thm3}
For Markovian rewards satisfying {\rm (C1)--(C3)},
SSMC achieves $E n_k^r = O(\log r)$ for $k \not\in \Xi$,
hence $R_N = O(\log N)$.
\end{thm}

Agrawal, Tenekatzis and Anantharam (1989) and Graves and Lai (1997) considered control problems in which,
instead of (\ref{PYB}) with $K$ Markov chains, 
there are $K$ arms with each arm representing a distinct Markov transition kernel acting on the same chain.
Tekin and Liu (2010) on the other hand considered (\ref{PYB}), 
with the constraints that $\cX$ is finite and $f_k(\cdot|x)$ is a point mass function for all $k$ and $x$.
They provided a UCB algorithm that achieves logarithmic regret.  

We can apply Theorem \ref{thm3} to show logarithmic regret for i.i.d. rewards on non-exponential parametric families.
Lai and Robbins (1985) showed that for the double exponential (DE) densities
\begin{equation} \label{DEfk}
f_k(y) = \tfrac{1}{2 \tau} e^{-|y-\mu_k|/\tau}, 
\end{equation}
with $\tau>0$,
efficiency is achieved by a UCB strategy involving KL-information of the DE densities,
hence implementation requires knowledge that the family is DE,
including knowing $\tau$.
In Example 1 below we state logarithmic regret, 
rather than efficiency, 
for SSMC.
The advantage of SSMC is that we do not assume knowledge of (\ref{DEfk}) in its implementation.
Verifications of (C1)--(C3) under (\ref{DEfk}) is given in Appendix B.

\smallskip
{\sc Example} 1.
For the double exponential densities {\rm (\ref{DEfk})}, 
conditions (C1)--(C3) hold,
hence under SSMC, 
$E n_k^r = O(\log r)$ for $k \not\in \Xi$.

\section{Numerical studies}

We compare SSMC and SSTC against procedures described in Section 2.1,
as well as more modern procedures like BESA,
KL-UCB,
UCB-Bayes and Thompson sampling.
The reader can refer to Chapters 1--3 of Kaufmann (2014) for a description of these procedures.
In Examples~2 and 3 we consider normal rewards and the comparisons are against procedures in which either efficiency or logarithmic regret has been established.
In Example 4 we consider double exponential rewards and there the comparisons are against procedures that have been shown to perform well numerically.
In Examples~5--7 we perform comparisons under the settings of Baransi, Maillard and Mannor (2014).

In the simulations done here $J=10000$ datasets are generated for each $N$, 
and the regret of a procedure is estimated by averaging over $\sum_{k=1}^K (\mu_*-\mu_k) N_k$.
Standard errors are located after the $\pm$ sign.
In Examples~5--7 we reproduce simulation results from Baransi et al. (2014).
Though no standard errors are provided,
they are likely to be small given that a larger $J=50000$ number of datasets are generated there.

\begin{table}[t]
\begin{center}
\begin{tabular}{c|rr}
& \multicolumn{2}{c}{Regret} \cr \cline{2-3}
& $N=1000$ & $N=10000$ \cr \hline
SSMC & 88.4$\pm$0.2 & 137.0$\pm$0.5 \cr
UCB1 & 90.2$\pm$0.3 & 154.4$\pm$0.7 \cr
UCB-Agrawal & 113.0$\pm$0.3 & 195.7$\pm$0.8 \cr
\end{tabular}
\caption{ The regrets of SSMC, UCB1 and UCB-Agrawal. 
The rewards have normal distributions with unit variances.
For each $N$ we generate $\mu_k \sim N(0,1)$ for $1 \leq k \leq 10$ a total of $J=10000$ times.}
\label{tab3}
\end{center}
\end{table}

\smallskip
{\sc Example} 2.
Consider $Y_{kt} \sim$ N($\mu_k,1$),
$1 \leq k \leq 10$.
In Table 1 we see that SSMC improves upon UCB1 and outperforms UCB-Agrawal [setting $b_n = \log \log \log n$ in (\ref{Agr})].
Here we generate $\mu_k \sim$ N(0,1) in each dataset. 

\begin{table}[t]
\begin{center}
\begin{tabular}{r|ll}
& \multicolumn{2}{c}{Regret} \cr \cline{2-3}
& $N=1000$ & $N=10000$ \cr \hline
SSTC & \ 239$\pm$1 & \ 492$\pm$5 \cr
UCB1-tuned & \ 130$\pm$2 & \ 847$\pm$23 \cr
UCB1-Normal & 1536$\pm$5 & 4911$\pm$31 
\end{tabular}
\caption{ The regrets of SSTC, UCB1-tuned and UCB1-Normal.
The rewards have normal distributions with unequal and unknown variances.
For each $N$ we generate $\mu_k \sim$ N(0,1) and $\sigma_k^{-2} \sim$ Exp(1) for $1 \leq k \leq 10$ a total of $J=10000$ times.}
\label{tab4}
\end{center}
\end{table}

\smallskip
{\sc Example} 3.
Consider $Y_{kt} \sim$ N($\mu_k,\sigma_k^2$), 
$1 \leq k \leq 10$.
We compare SSTC against UCB1-tuned and UCB1-Normal.
UCB1-tuned was suggested by Auer et al. and shown to perform well numerically.
Under UCB1-tuned the population $\Pi_k$ maximizing
$$\bar Y_{k n_k} + \sqrt{\tfrac{\log n}{n_k} \min(\tfrac{1}{4}, V_{kn})},
$$
where $V_{kn} = \wht \sigma_{k n_k}^2+\sqrt{\frac{2 \log n}{n_k}}$, 
is selected.
In Table 2 we see that UCB1-tuned is significantly better at $N=1000$ whereas SSTC is better at $N=10000$. 
UCB1-Normal performs quite poorly. 
Here we generate $\mu_k \sim {\rm N}(0,1)$ and $\sigma_k^{-2} \sim {\rm Exp}(1)$ in each dataset.

Kaufmann, Capp\`{e} and Garivier (2012) performed simulations under the setting of normal rewards with unequal variances,
with $(\mu_1,\sigma_1)=(1.8,0.5)$,
$(\mu_2,\sigma_2)=(2,0.7)$,
$(\mu_3,\sigma_3)=(1.5,0.5)$ and $(\mu_4,\sigma_4)=(2.2,0.3)$.
They showed that UCB-Bayes achieves regret of about 28 at $N=1000$ and about 47 at $N=10000$.
We apply SSTC on this setting,
achieving regrets of 26.0$\pm$0.1 at $N=1000$ and 43.3$\pm$0.2 at $N=10000$.

\begin{table}[t]
\begin{center}
\begin{tabular}{r|lll|lll}
& \multicolumn{3}{c|}{Regret} & \multicolumn{3}{c}{Regret $(\times 10)$} \cr \cline{2-7}
& \multicolumn{3}{c|}{$N=1000$} & \multicolumn{3}{c}{$N=10000$} \cr \cline{2-7}
& $\lambda=1$ & $\lambda=2$ & $\lambda=5$ & $\lambda=1$ & $\lambda=2$ & $\lambda=5$ \cr \hline
SSMC & 141.7$\pm$0.4 & 330$\pm$1 & 795$\pm$3 & \ 23.6$\pm$0.1 & \ 65.0$\pm$0.3 & 236.9$\pm$0.8 \cr 
BESA &  117$\pm$1 & 265$\pm$2 & 627$\pm$3 & \ 28.9$\pm$0.7 & \ 73$\pm$1 & 215$\pm$2 \cr
UCB1-tuned &  101$\pm$2 & 244$\pm$3 & 608$\pm$6 & \ 50$\pm$1 & 183$\pm$3 & 499$\pm$6 \cr \hline
Boltz $\tau=$0.1 & 130$\pm$2 & 294$\pm$4 & 673$\pm$7 & \ 84$\pm$2 & 224$\pm$4 & 557$\pm$6 \cr
0.2 & 128$\pm$2 & 264$\pm$3 & 632$\pm$6 & \ 80$\pm$1 & 169$\pm$3 & 465$\pm$6 \cr
0.5 & 332$\pm$1 & 387$\pm$2 & 632$\pm$5 & 310$\pm$5 & 311$\pm$2 & 428$\pm$4 \cr
1 & 728$\pm$2 & 737$\pm$2 & 816$\pm$4 & 731$\pm$2 & 716$\pm$2 & 712$\pm$3 \cr \hline
$\epsilon$-greedy $c=$0.1 & 170$\pm$3 & 327$\pm$4 & 681$\pm$7 & 133$\pm$3 & 283$\pm$4 & 579$\pm$7 \cr
0.2 & 162$\pm$3 & 312$\pm$4 & 653$\pm$6 & 114$\pm$2 & 251$\pm$4 & 536$\pm$6 \cr
0.5 & 150$\pm$2 & 282$\pm$3 & 604$\pm$6 & \ 82$\pm$2 & 189$\pm$3 & 444$\pm$5 \cr
1 & 159$\pm$2 & 271$\pm$3 & 569$\pm$5 & \ 61$\pm$1 & 146$\pm$3 & 370$\pm$5 \cr
2 & 200$\pm$1 & 289$\pm$2 & 559$\pm$4 & \ 52.9$\pm$0.9 & 113$\pm$2 & 302$\pm$4 \cr
5 & 334$\pm$1 & 396$\pm$2 & 617$\pm$4 & \ 63.4$\pm$0.5 & 101$\pm$1 & 241$\pm$3 \cr
10 & 524$\pm$2 & 567$\pm$2 & 742$\pm$3 & \ 95.7$\pm$0.4 & 119.5$\pm$0.8 & 226$\pm$2 \cr
20 & 811$\pm$3 & 839$\pm$3 & 951$\pm$3 & 156.9$\pm$0.5 & 172.1$\pm$0.7 & 251$\pm$2 
\end{tabular}
\end{center}
\caption{ Regret comparisons for double exponential density rewards.
For each $N$ and $\lambda$ we generate $\mu_k \sim$ N(0,1) for $1 \leq k \leq 10$ a total of $J=10000$ times.}
\end{table}

\smallskip
{\sc Example} 4.
Consider double exponential rewards $Y_{kt} \sim f_k$,
with densities
$$f_k(y) = \tfrac{1}{2 \lambda} e^{-|y-\mu_k|/\lambda}, \quad 1 \leq k \leq 10.
$$
We compare SSMC against UCB1-tuned, 
BESA,
Boltzmann exploration and $\epsilon$-greedy.
For $\epsilon$-greedy we consider $\epsilon_n = \min(1,\tfrac{3c}{n})$.
We generate $\mu_k \sim$ N(0,1) in each dataset.

Table 3 shows that UCB1-tuned has the best performances at $N=1000$, 
whereas SSMC has the best performances at $N=10000$.
BESA does well for $\lambda=2$ at $N=1000$,
and also for $\lambda=5$ at $N=10000$.
A properly-tuned Boltzmann exploration does well at $N=1000$ for $\lambda=2$,
whereas a properly-tuned $\epsilon$-greedy does well at $\lambda=2$ and 5 for $N=1000$ and at $\lambda=5$ for $N=10000$.

\begin{table}[t]
\begin{center}
\begin{tabular}{c|r|r|r|r|r|r|r|r}
& \multicolumn{7}{c|}{Frequency of emp. regrets} \cr 
& \multicolumn{7}{c|}{lying within a given range} \cr \cline{2-8}
& 0 & 200 & 400 & 600 & 800 & 1000 & 1200 \cr
& to & to & to & to & to & to & to & Worst \cr
& 200 & 400 & 600 & 800 & 1000 & 1200 & 2100 & emp. regret\cr \hline
SSMC & 9134 & 845 & 16 & 5 & 0 & 0 & 0 &  770 \cr
BESA & 9314 & 424 & 143 & 66 & 27 & 15 & 11 & 2089 \cr
UCB1-tuned &  8830 & 625 & 301 & 132 & 64 & 32 & 16 & 1772    
\end{tabular}
\end{center}
\caption{ Number of simulations $($out of $10000)$ lying within a given empirical regret range,
and the worst empirical regret, 
when $N=1000$ and $\lambda=1$.}
\end{table}

\begin{table}[t]
\begin{center}
\begin{tabular}{c|r|r|r|r|r|r|r|r}
& \multicolumn{7}{c|}{Frequency of emp. regrets} \cr 
& \multicolumn{7}{c|}{lying within a given range} \cr \cline{2-8}
& 0 & 1000 & 2000 & 3000 & 4000 & 5000 & 10000 & \cr
& to & to & to & to & to & to & to & Worst \cr
& 1000 & 2000 & 3000 & 4000 & 5000 & 10000 & 21000 & emp. regret\cr \hline
SSMC & 9988 & 8 & 3 & 0 & 0 & 1 & 0 & 6192 \cr
BESA & 9708 & 125 & 59 & 34 & 25 & 40 & 9 & 20639 \cr
UCB1-tuned & 8833 & 365 & 250 & 161 & 122 & 225 & 44 & 16495    
\end{tabular}
\end{center}
\caption{ Number of simulations $($out of $10000)$ lying within a given empirical regret range,
and the worst empirical regret, 
when $N=10000$ and $\lambda=1$.}
\end{table}

In Tables 4 and 5 we tabulate the frequencies of the empirical regrets $\sum_{k=1}^K (\mu_*-\mu_k) N_k$ over the $J=10000$ simulation runs each for $N=1000$ and 10000, 
at $\lambda=1$,
for SSMC,
BESA and UCB1-tuned.
Tha tables show that SSMC has the best control of excessive sampling of inferior arms, 
the worst empirical regret being less than half that of BESA and UCB1-tuned.
 
\begin{table}[t]
\begin{center}
\begin{tabular}{l|llll}
& \multicolumn{4}{c}{Scenario} \cr \cline{2-5}
& 1 & 2 & 3 & 4 \cr \hline
SSMC & 12.4$\pm$0.1 & 43.1$\pm$0.4 & \ 97.9$\pm$0.2 & 165.3$\pm$0.2 \cr
SSMC$^*$ & \ 9.5$\pm$0.2 & 48.5$\pm$0.6 & \ 64.4$\pm$0.3 & 156.0$\pm$0.4 \cr \hline
BESA & 11.83 & 42.6 & \ 74.41 & 156.7 \cr
KL-UCB & 17.48 & 52.34 & 121.21 & 170.82 \cr
KL-UCB$+$ & 11.54 & 41.71 & \ 72.84 & 165.28 \cr
Thompson & 11.3 & 46.14 & \ 83.36 & 165.08 \cr
\end{tabular}
\end{center}
\caption{ Regret comparisons for Bernoulli rewards.}
\end{table}

\smallskip
{\sc Example} 5.
Consider $N=20000$ Bernoulli rewards under the following scenarios. 

\begin{enumerate}
\item $\mu_1=0.9$, $\mu_2=0.8$.
\item $\mu_1=0.81$, $\mu_2=0.8$.
\item $\mu_2=0.1$, $\mu_2=\mu_3=\mu_4=0.05$, 
$\mu_5=\mu_6=\mu_7=0.02$, 

$\mu_8=\mu_9=\mu_{10}=0.01$.
\item $\mu_1=0.51$, $\mu_2=\cdots=\mu_{10}=0.5$.
\end{enumerate}

\noindent When comparing the simulated regrets in Table 6,
it is useful to remember that BESA and SSMC are non-parametric,
using the same procedures even when the rewards are not Bernoulli,
whereas KL-UCB and Thompson sampling utilize information on the Bernoulli family.
SSMC$^*$ is a variant of SSMC,
see Section~6,
with more moderate levels of explorations. 

\begin{table}[t]
\begin{center}
\begin{tabular}{l|ll}
& Trunc. expo. & Trunc. Poisson \cr \hline
SSMC & 33.8$\pm$0.4 & 18.6$\pm$0.1 \cr
SSMC$^*$ & 29.6$\pm$0.7 & 14.7$\pm$0.2 \cr \hline
BESA & 53.26 & 19.37 \cr
BESAT & 31.41 & 16.72 \cr
KL-UCB-expo & 65.67 & --- \cr
KL-UCB-Poisson & --- & 25.05 \cr
\end{tabular}
\end{center}
\caption{ Regret comparisons for truncated exponential and Poisson rewards.}
\end{table}

\smallskip
{\sc Example} 6. 
Consider truncated exponential and Poisson distributions with $N=20000$. 
For truncated exponential we consider $Y_{kt} = \min(\frac{X_{kt}}{10},1)$,
where $X_{kt} \stackrel{\rm i.i.d.}{\sim} {\rm Exp}(\lambda_k)$ (density $\lambda_k e^{-\lambda_k x}$)
with $\lambda_k = \frac{1}{k}$, $1 \leq k \leq 5$.
For truncated Poisson we consider $Y_{kt} = \min(\frac{X_{kt}}{10},1)$,
where $X_{kt} \stackrel{\rm i.i.d.}{\sim} {\rm Poisson}(\lambda_k)$,
with $\lambda_k = 0.5+\frac{k}{3}$, $1 \leq k \leq 6$.
The simulation results are given in Table 7.
BESAT is a variation of BESA that starts with 10 observations from each population.

\smallskip
{\sc Example} 7. 
Consider $K=2$ and $N=20000$ with $Y_{1t} \stackrel{\rm i.i.d.}{\sim} {\rm Uniform}(0.2,0.4)$ and $Y_{2t} \stackrel{\rm i.i.d.}{\sim} {\rm Uniform}(0,1)$.
Here SSMC underperforms with regret of 163$\pm$7 compared to Thompson sampling,
which has regret of 13.18.
On the other hand SSTC,
by normalizing the different scales of the two uniform distributions,
is able to achieve the best regret of 2.9$\pm$0.2.

\section{Discussion}

Together with BESA,
the procedures SSMC and SSTC that we introduce here form a class of non-parametric procedures that differ from traditional non-parametric procedures,
like $\epsilon$-greedy and Boltzmann exploration, 
in their recognition that when deciding between which of two populations to be sampled,
samples or subsamples of the same rather than different sizes should be compared.
Among the parametric procedures, 
Thompson sampling fits most with this scheme.

As mentioned earlier, 
in SSMC (and SSTC),
when the leading population $\Pi_{\zeta}$ in the previous round is sampled,
essentially only one additional comparison is required in the current round between $\Pi_{\zeta}$ and $\Pi_k$ for $k \neq \zeta$.
On the other hand when there are $n$ rewards, 
an order $n$ comparisons may be required between $\Pi_{\zeta}$ and $\Pi_k$ when $\Pi_k$ wins in the previous round.
It is these added comparisons that,
relative to BESA,
allows for faster catching-up of a potentially undersampled optimal arm.
Tables~4 and 5 show the benefits of such added explorations in minimizing the worst-case empirical regret.

To see if SSMC still works well if we moderate these added explorations,  
we experimented with the following variation of SSMC in Examples 6 and~7. 
The numerical results indicate improvements.

\smallskip
\underline{SSMC$^*$}

Proceed as in SSMC, 
with step 2(b)iii. replaced by the following.

\smallskip
2(b)iii$'$ If $c_n \leq n_k < n_{\zeta}$,
then $\Pi_k$ wins the challenge when
$$\bar Y_{k n_k} \geq \bar Y_{\zeta,t:(t+n_k-1)} \mbox{ for some } t = 1+un_k, 0 \leq u \leq \lfloor \tfrac{n_{\zeta}}{n_k} \rfloor-1.
$$

\smallskip
In contrast to SSMC,
in SSMC$^*$ we partition the rewards of the leading arm into groups of size $n_k$ for comparisons instead of reusing the rewards in moving-averages.
In principle the members of the group need not be consecutive in time,
thus allowing for the modifications of SSMC$^*$ to provide storage space savings when the support of the distributions is finite.
That is rather than to store the full sequence,
we simply store the number of occurrences at each support point,
and generate a new (permuted) sequence for comparisons whenever necessary. 
Likewise in BESA,
there is substantial storage space savings for finite-support distributions by storing the number of occurrences at each support point.

\section{Proofs of Theorems \ref{thm1}--\ref{thm3}}

Since SSMC and SSTC are index-blind,
we may assume without loss of generality that $\mu_1 = \mu_*$.
We provide here the statements and proofs of supporting Lemmas \ref{lem2} and \ref{lem1},
and follow up with the proofs of Theorems~\ref{thm1}--\ref{thm3} in Sections~7.1--7.3.
We denote the complement of an event $D$ by $\bar D$,
let $\lfloor \cdot \rfloor$ and $\lceil \cdot \rceil$ denote the greatest and least integer function respectively,
and let $|A|$ denote the number of elements in a set $A$.

Let $n_k^r(=n_k)$ be the number of observations from $\Pi_k$ at the beginning of round $r$.
Let $n^r(=n)= \sum_{k=1}^K n_k^r$. 
Let $n_*^r = \max_{1 \leq k \leq K} n_k^r$.
Let 
\begin{eqnarray*}
& & \Xi = \{ \ell: \mu_{\ell} = \mu_* \} \mbox{ be the set of optimal arms}, \cr
& & \zeta^r(=\zeta) \mbox{ the leader at the beginning of round } r(\geq 2).
\end{eqnarray*}
More specifically,
let
\begin{eqnarray*}
\cZ^r & = & \{ k: n_k^r = n_*^r \}, \cr
\cZ_1^r & = & \{ \ell \in \cZ^r: \bar Y_{\ell n_{\ell}^r} \geq \bar Y_{k n_k^r} \mbox{ for all } k \in \cZ^r \}.
\end{eqnarray*}
If $\zeta^{r-1} \in \cZ_1^r$,
then $\zeta^r=\zeta^{r-1}$.
Otherwise the leader $\zeta^r$ is selected randomly (uniformly) from $\cZ_1^r$.
In particular if $\cZ_1^r$ has a single element,
then that element must be $\zeta^r$.
For $r \geq 2$,
let
$$A^r = \{ \zeta^r \not\in \Xi \} = \{ \mbox{leader at round } r \mbox{ is inferior} \}.
$$
We restrict to $r \geq 2$ because the leader is not defined at $r=1$.
Likewise in our subsequent notations on events $B^r$, $C^r$, $D^r$, $G_k^r$ and $H_k^r$,
we restrict to $r \geq 2$.

In Lemma \ref{lem2} below the key ingredient leading to (\ref{xik}) is condition (I) on the event $G_k^r$,
which says that it is difficult for an inferior arm $k$ with at least $(1+\epsilon) \xi_k \log r$ rewards to win against a leading optimal arm $\zeta$.
In the case of exponential families we show efficiency by verifying (I) with $\xi_k = \frac{1}{I_1(\mu_k)}$.
Condition (II),
on the event $H_k^r$,
says that analogous winnings from an inferior arm $k$ with at least $J_k \log r$ rewards,
for $J_k$ large,
are asymptotically negligible. 
Condition  (III) limits the number of times an inferior arm is leading.
This condition is important because $G_k^r$ and $H_k^r$ refer to the winning of arm $k$ when the leader is optimal,
hence the need, 
in (III), 
to bound the event probability of an inferior leader.

\begin{lem} \label{lem2}
Let $k \not\in \Xi$ $($i.e. $k$ is not an optimal arm$)$ and define
\begin{eqnarray} \label{Gkr}
G_k^r & = & \{ \zeta^s \in \Xi, n_k^{s+1}=n_k^s+1, \\ \nonumber
& & \quad n_k^s \geq (1+\epsilon) \xi_k \log r \mbox{ for some } 2 \leq s \leq r-1 \}, \\ \label{Hkr}
H_k^r & = & \{ \zeta^s \in \Xi, n_k^{s+1}=n_k^s+1, \\ \nonumber
& & \quad n_k^s \geq J_k \log r \mbox{ for some } 2 \leq s \leq r-1 \},
\end{eqnarray}
 for some $\epsilon>0$,
$\xi_k>0$ and $J_k>0$.
Consider the following conditions.

\smallskip
{\rm (I)} There exists $\xi_k>0$ such that for all $\epsilon>0$,
$P(G_k^r) \rightarrow 0$ as $r \rightarrow \infty$.

\smallskip
{\rm (II)} There exists $J_k>0$ such that $P(H_k^r) = O(r^{-1})$ as $r \rightarrow \infty$.

\smallskip
{\rm (III)} $P(A^r) = o(r^{-1})$ as $r \rightarrow \infty$.

\smallskip
\noindent Under {\rm (I)--(III)},
\begin{equation} \label{xik}
\limsup_{r \rightarrow \infty} \frac{En_k^r}{\log r} \leq \xi_k.
\end{equation}
\end{lem}

\medskip
{\sc Proof}.
Consider $r \geq 3$. 
Let $b_r = 1+(1+\epsilon) \xi_k \log r$ and $d_r = 1+J_k \log r$.
Under the event $\bar G_k^r$,
arm $k$ in round $s \in [2,r-1]$ is sampled to a size beyond $b_r$ only when $\zeta^s \not\in \Xi$ (i.e. under the event $A^s$).
In view that $n_k^2=1(<b_r)$,
it follows that
$$n_k^r \leq b_r + \sum_{s=2}^{r-1} {\bf 1}_{A^s}.
$$ 
Hence
\begin{equation} \label{nG}
n_k^r {\bf 1}_{\bar G_k^r} \leq b_r + \sum_{s=2}^{r-1} {\bf 1}_{A^s}.
\end{equation}
Similarly under the event $\bar H_k^r$, 
$$n_k^r \leq d_r+\sum_{s=2}^{r-1} {\bf 1}_{A^s}.
$$ 
Hence
\begin{equation} \label{nH}
n_k^r {\bf 1}_{(G_k^r \setminus H_k^r)} \leq d_r {\bf 1}_{G_k^r} + \sum_{s=2}^{r-1} {\bf 1}_{A^s}.
\end{equation}
Since $n_k^r \leq r$,
by (\ref{nG}) and (\ref{nH}),
\begin{eqnarray} \label{Enkr}
En_k^r & = & E(n_k^r {\bf 1}_{G_k^r \cap H_k^r}) + E(n_k^r {\bf 1}_{(G_k^r \setminus H_k^r)}) + E(n_k^r {\bf 1}_{\bar G_k^r}) \\ \nonumber
& \leq & r P(H_k^r) + \Big[ d_r P(G_k^r)+\sum_{s=2}^{r-1} P(A^s) \Big] + \Big[ b_r+\sum_{s=2}^{r-1} P(A^s) \Big].
\end{eqnarray} 
By (III),
$\sum_{s=2}^r P(A^s) = o(\log r)$,
therefore by (\ref{Enkr}), 
(I) and (II), 
$$\limsup_{r \rightarrow \infty} \frac{En_k^r}{\log r} \leq (1+\epsilon) \xi_k.
$$
We can thus conclude (\ref{xik}) by letting $\epsilon \rightarrow 0$.
$\wbox$  

\medskip
The verification of (III) is made easier by Lemma \ref{lem1} below.
To provide intuitions for the reader we sketch its proof first before providing the details.

\begin{lem} \label{lem1}
Let 
\begin{eqnarray*}
B^s & = & \{ \zeta^s \in \Xi, n_k^{s+1} = n_k^s+1, 
n_k^s = n_{\zeta}^s-1 \mbox{ for some } k \not\in \Xi \}, \cr
C^s & = & \{ \zeta^s \not\in \Xi, n_{\ell}^{s+1} = n_{\ell}^s \mbox{ for some } \ell \in \Xi \}.
\end{eqnarray*}
If as $s \rightarrow \infty$,
\begin{eqnarray} \label{BC}
P(B^s) & = & o(s^{-2}), \\ \label{Cs}
P(C^s) & = & o(s^{-1}),
\end{eqnarray}
then $P(A^r) = o(r^{-1})$ as $r \rightarrow \infty$.
\end{lem}

{\sc Sketch of proof}. 
Note that (\ref{BC}) bounds the probability of an inferior arm taking the leadership from an optimal leader in round $s+1$,
whereas (\ref{Cs}) bounds the probability of an inferior leader winning against an optimal challenger in round $s$.
Let $s_0 = \lfloor \frac{r}{4} \rfloor$ and for $r \geq 8$,
let 
\begin{eqnarray*}
D^r & = & \{ \zeta^s \in \Xi \mbox{ for some } s_0 \leq s \leq r-1 \} \cr
& = & \{ \mbox{the leader is optimal for some rounds between } s_0 \mbox{ to } r-1 \}.
\end{eqnarray*}

Under $A^r \cap D^r$,
there is a leadership takeover by an inferior arm at least once between rounds $s_0+1$ and $r$.
More specifically let $s_1$ be the largest $s \in [s_0,r-1]$ for which $\zeta^s \in \Xi$.
If $s_1 < r-1$,
then by the definition of $s_1$,
$\zeta^{s_1+1} \not\in \Xi$.
If $s_1=r-1$,
then since we are under $A^r$,
$\zeta^{s_1+1}=\zeta^r \not\in \Xi$.
In summary
\begin{eqnarray} \label{ArDr}
A^r \cap D^r & = & \{ \xi^s \in \Xi \mbox{ for some } s_0 \leq s \leq r-1, \zeta^r \not\in \Xi \} \\ \nonumber
& \subset & \cup_{s=s_0}^{r-1}  \{ \zeta^s \in \Xi, \zeta^{s+1} \not\in \Xi \}.
\end{eqnarray}
By showing that
\begin{equation} \label{show}
\{ \zeta^s \in \Xi, \zeta^{s+1} \not\in \Xi \} \subset B^s,
\end{equation}
we can conclude from (\ref{BC}) and (\ref{ArDr}) that 
\begin{equation} \label{AD}
P(A^r \cap D^r) \leq \sum_{s=s_0}^{r-1} P(B^s) = o(rs_0^{-2}) = o(r^{-1}).
\end{equation}
To see (\ref{show}),
recall that by step 2(b)i of SSMC or SSTC, 
if the (optimal) leader and (inferior) challenger have the same sample size,
then the challenger loses by default.
The tie-breaking rule then ensures that the challenger is unable to take over leadership in the next round.
Hence for $\zeta^s$ to lose leadership to an inferior arm $k$ in round $s+1$,
it has to lose to arm $k$ when arm $k$ has exactly $n_{\zeta}^s-1$ observations.

What (\ref{AD}) says is that if at some previous round $s \geq s_0$ the leader is optimal,
then (\ref{BC}) makes it difficult for an inferior arm to take over leadership during and after round $s$,
so the leader is likely to be optimal all the way from rounds $s$ to $r$.
The only situation we need to guard against is $\bar D^r$,
the event that leaders are inferior for all rounds between $s_0$ and $r-1$.
Let $\#^r = \sum_{s=s_0}^{r-1} {\bf 1}_{C^s}$ be the number of rounds an inferior leader wins against at least one optimal arm.
In (\ref{Hr}) we show that by (\ref{Cs}),
the optimal arms will,
with high probability,
lose less than $\frac{r}{4}$ times between rounds $s_0$ and $r-1$ when the leader is inferior.

We next show that 
\begin{equation} \label{Drc}
\bar D^r \subset \{ \#^r \geq \tfrac{r}{4} \},
\end{equation}
(or $\{ \#^r < \frac{r}{4} \} \subset D^r$), 
that is if the optimal arms lose this few times,
then one of them has to be a leader at some round between $s_0$ to $r-1$.
Lemma~\ref{lem1} follows from (\ref{AD})--(\ref{Hr}).

\smallskip
{\sc Proof of Lemma} \ref{lem1}. 
Consider $r \geq 8$.
By (\ref{Cs}),
$$E(\#^r) = \sum_{s=s_0}^{r-1} P(C^s) = o(rs_0^{-1}) \rightarrow 0,
$$
hence by Markov's inequality,
\begin{equation} \label{Hr}
P(\#^r \geq \tfrac{r}{4}) \leq \tfrac{E(\#^r)}{r/4} = o(r^{-1}).
\end{equation}

It remains for us to show (\ref{Drc}).
Assume $\bar D^r$.
Let $m^s = n_{\zeta}^s-\max_{\ell \in \Xi} n_{\ell}^s$.
Observe that $n_{\zeta}^{s+1}=n_{\zeta}^s$ if $n_{\ell}^{s+1}=n_{\ell}^s+1$ for some $\ell \neq \zeta^s$. 
This is because the leader $\zeta^s$ is not sampled if it loses at least one challenge.
Moreover by step 2(b)i. of SSMC or SSTC, 
all arms with the same number of observations as $\zeta^s$ are not sampled.
Therefore if $\zeta^s \not\in \Xi$ and $n_{\ell}^{s+1}=n_{\ell}^s+1$ for all $\ell \in \Xi$,
that is if all optimal arms win against an inferior leader,
then $m^{s+1} = m^s-1$.
In other words,
\begin{equation} \label{Fs}
F^s := \{ \zeta^s \not\in \Xi, n_{\ell}^{s+1}=n_{\ell}^s+1 \mbox{ for all } \ell \in \Xi \} \subset \{ m^{s+1}=m^s-1 \}.
\end{equation}
Since $m^{s+1} \leq m^s+1$,
it follows from (\ref{Fs}) that
$m^{s+1} \leq m^s+1-2 {\bf 1}_{F^s}$.
Therefore 
$$m^r \leq m^{s_0}+(r-s_0)-2 \sum_{s=s_0}^{r-1} {\bf 1}_{F^s}, \\ \nonumber
$$
and since $m^r \geq 0$ and $m^{s_0} \leq s_0$,
we can conclude that
\begin{equation} \label{Fbound}
\sum_{s=s_0}^{r-1} {\bf 1}_{F^s} \leq \tfrac{r}{2}.
\end{equation}

Under $\bar D^r$,
${\bf 1}_{C^s} = 1-{\bf 1}_{F^s}$ for $s_0 \leq s \leq r-1$,
and it follows from (\ref{Fbound}) that 
$$\#^r \geq (r-s_0)-\tfrac{r}{2} \geq \tfrac{r}{4},
$$
and (\ref{Drc}) indeed holds.
$\wbox$

\subsection{Proof of Theorem \ref{thm1}}
We consider here SSMC. 
Equation (\ref{BC}) follows from Lemma \ref{lem3a} below and $c_r=o(\log r)$ whereas (\ref{Cs}) follows from Lemma \ref{lem3} and $\frac{c_r}{\log \log r} \rightarrow \infty$.
We can thus conclude $P(A^r)=o(r^{-1})$ from Lemma \ref{lem1},
and together with the verification in Lemma \ref{lem4} of (I),
see Lemma \ref{lem2},
for $\xi_k=1/I_1(\mu_k)$ and (II) for $J_k$ large, 
we can conclude Theorem \ref{thm1}.

The proofs of Lemmas \ref{lem3a}--\ref{lem4} use large deviations Chernoff bounds that are given below in Lemma \ref{lem3b}.
They can be shown using change-of-measure arguments.
Let $I_k$ be the large deviations rate function of $f_k$.

\begin{lem} \label{lem3b}
Under {\rm (\ref{expo})},
if $1 \leq k \leq K$,
$t \geq 1$ and $\omega = \psi'(\theta)$ for some $\theta \in \Theta$,
then
\begin{eqnarray} \label{LD1}
P(\bar Y_{kt} \geq \omega) & \leq & e^{-t I_k(\omega)} \mbox{ if } \omega > \mu_k, \\
\label{LD2}
P(\bar Y_{kt} \leq \omega) & \leq & e^{-t I_k(\omega)} \mbox{ if } \omega < \mu_k.
\end{eqnarray}
\end{lem}

In Lemmas \ref{lem3a}--\ref{lem4} we let $\omega=\frac{1}{2}(\mu_*+\max_{k:\mu_k < \mu_*} \mu_k)$ and $a=\min_{1 \leq k \leq K} I_k(\omega)$.
Recall that the parameter $c_r$ is a threshold for forced explorations,
in step 2(b)ii. of SSMC. 

\begin{lem} \label{lem3a}
Under {\rm (\ref{expo})},
$P(B^r) \leq \frac{3K^2}{1-e^{-a}} e^{-a(\frac{r}{K}-1)}$ when $\frac{r}{K}-1 \geq c_r$. 
\end{lem}

{\sc Proof}.
Let $r$ be such that $\frac{r}{K}-1 \geq c_r$.
The event $B^r$ occurs if at round $r$ the leading arm $\ell$ is optimal (i.e. $\ell \in \Xi$), 
and it loses to an inferior arm $k (\not\in \Xi)$ with $n_k =u$ and $n_{\ell}=u+1$ for $u+1 \geq \frac{r}{K}$ (since arm $\ell$ is leading). 
It follows from Lemma \ref{lem3b} that
\begin{eqnarray*} 
P(\bar Y_{\ell,t:(t+u-1)} \leq \omega \mbox{ for } t =1 \mbox{ or } 2) & \leq & 2 e^{-u I_{\ell}(\omega)}, \quad \ell \in \Xi, \cr
P(\bar Y_{k u} \geq \omega) & \leq & e^{-u I_k(\omega)}, \quad k \not\in \Xi.
\end{eqnarray*}
Since arm $\ell$ loses to arm $k$ when $\bar Y_{ku} \geq \min(\bar Y_{\ell,1:u}, \bar Y_{\ell,2:(u+1)})$,
it follows that
$$P(B^r) \leq \sum_{\ell \in \Xi} \sum_{k \not\in \Xi} \sum_{u = \lceil \frac{r}{K} \rceil-1}^r (2e^{-u I_\ell(\omega)}+e^{-u I_k(\omega)}), 
$$
and Lemma \ref{lem3a} holds. 
$\wbox$

\begin{lem} \label{lem3}
Under {\rm (\ref{expo})}, 
$P(C^r) \leq K^2 e^{-c_r a} \frac{(\log r)^6}{r} + o(r^{-1})$. 
\end{lem}

{\sc Proof}. 
The event $C^r$ occurs if at round $r$ the leading arm $k$ is inferior (i.e. $k \not\in \Xi$), 
and it wins a challenge against one or more optimal arms $\ell (\in \Xi)$. 
By step 2(b)ii. of SSMC,
arm~$k$ loses automatically when $n_{\ell} < c_n$,
hence we need only consider $n_{\ell} \geq c_n$.
Note that when $n_k=n_{\ell}$,
for arm $k$ to be the leader,
by the tie-breaking rule we require $\bar Y_{k n_{\ell}} \geq \bar Y_{\ell n_{\ell}}$.
We shall consider $n_{\ell} > (\log r)^2$ in case 1 and $n_{\ell} = v$ for $c_n \leq v < (\log r)^2$ in case 2.

Case 1: $n_{\ell} > (\log r)^2$.
By Lemma \ref{lem3b},
\begin{eqnarray} \label{lem4.1}
P(\bar Y_{\ell n_{\ell}} \leq \omega \mbox{ for some } n_{\ell} > (\log r)^2) & \leq & \tfrac{1}{1-e^{-a}} e^{-a(\log r)^2} \\ \label{lem4.2}
P(\bar Y_{kn_{\ell}} \geq \omega \mbox{ for some } n_{\ell} > (\log r)^2) & \leq & \tfrac{1}{1-e^{-a}} e^{-a(\log r)^2}.
\end{eqnarray} 

Case 2: $n_{\ell}=v$ for $(c_r \leq) c_n \leq v < (\log r)^2$.
In view that $n_k \geq \frac{r}{K}$ when $k$ is the leading arm,
we shall show that for $r$ large,
for each such $v$ there exists $\xi (=\xi_v)$ such that
\begin{eqnarray} \label{eq1}
& & P(\bar Y_{\ell v} < \xi) \leq e^{-c_r a} \tfrac{(\log r)^4}{r}, \\ \label{eq2}
& & P(\bar Y_{k,t:(t+v-1)} > \xi \mbox{ for } 1 \leq t \leq \tfrac{r}{K}) \\ \nonumber
& [\leq & P(\bar Y_{k v} > \xi)^{\lfloor \frac{r}{K v} \rfloor}] \leq \exp[-\tfrac{(\log r)^2}{K}+1]. 
\end{eqnarray}
The inequality within the brackets in (\ref{eq2}) follows from partitioning $[1,\frac{r}{K}]$ into $\lfloor \frac{r}{K v} \rfloor$ segments of length $v$,
and applying independence of the sample on each segment.

Since $\theta_{\ell} > \theta_k$,
if $\sum_{t=1}^v y_t \leq v \mu_k$, 
then by (\ref{expo}),
\begin{eqnarray*}
\prod_{t=1}^v f(y_t;\theta_{\ell}) & = & e^{(\theta_{\ell}-\theta_k) \sum_{t=1}^v y_t - v[\psi(\theta_{\ell})-\psi(\theta_k)]} \prod_{t=1}^v f(y_t;\theta_k) \cr
& \leq & e^{-v I_{\ell}(\mu_k)} \prod_{t=1}^v f(y_t;\theta_k).
\end{eqnarray*}
Hence if $\xi \leq \mu_k$,
then as $v \geq c_r$,
\begin{equation} \label{eq3}
P(\bar Y_{\ell v} < \xi) \leq e^{-v I_{\ell}(\mu_k)} P(\bar Y_{k v} < \xi) \leq e^{-c_r a} P(\bar Y_{k v} < \xi).
\end{equation}
Let $\xi (\leq \mu_k$ for large $r$) be such that
\begin{equation} \label{eq4}
P(\bar Y_{k v} < \xi) \leq \tfrac{(\log r)^4}{r} \leq P(\bar Y_{k v} \leq \xi).
\end{equation}
Equation (\ref{eq1}) follows from (\ref{eq3}) and the first inequality in (\ref{eq4}),
whereas (\ref{eq2}) follows from the second inequality in (\ref{eq4}) and $v < (\log r)^2$. 
By (\ref{lem4.1})--(\ref{eq2}),
$$P(C^r) \leq \sum_{\ell \in \Xi} \sum_{k \not\in \Xi} \Big\{ \tfrac{2}{1-e^{-a}} e^{-a(\log r)^2} 
+ \sum_{v = \lceil c_r \rceil}^{\lfloor (\log r)^2 \rfloor} (e^{-c_r a} \tfrac{(\log r)^4}{r} + \exp\Big[-\tfrac{(\log r)^2}{K} +1 \Big] \Big) \Big\},
$$
and Lemma \ref{lem3} holds. 
$\wbox$

\begin{lem} \label{lem4}
Under {\rm (\ref{expo})} and $c_r=o(\log r)$,
{\rm (I)} $($in the statement of Lemma~{\rm {\ref{lem2}}}$)$ holds for $\xi_k = 1/I_1(\mu_k)$ and {\rm (II)} holds for $J_k > \max(\frac{1}{I_k(\omega)},\frac{2}{I_1(\omega)})$,
where $\omega = \frac{1}{2}(\mu_*+\max_{k:\mu_k < \mu_*} \mu_k)$.
\end{lem}

{\sc Proof}. 
Let $k \not\in \Xi$.
Let $\mu_k < \omega_k < \mu_1$ be such that $(1+\epsilon) I_1(\omega_k) > I_1(\mu_k)$.
Consider $n_k=u$ for $u \geq (1+\xi_k) \log r$ (in $G_k^r$) and $u \geq J_k \log r$ (in $H_k^r$).
Since $I_{\ell} = I_1$ for $\ell \in \Xi$,
it follows from Lemma \ref{lem3b} that 
\begin{eqnarray} \label{4.1}
P(\bar Y_{{\ell},t:(t+u-1)} \leq \omega_k \mbox{ for some } 1 \leq t \leq r) & \leq & re^{-u I_1(\omega_k)}, \\ \label{4.2} 
P(\bar Y_{ku} \geq \omega_k) & \leq & e^{-u I_k(\omega_k)}.
\end{eqnarray}

Since $c_r = o(\log r)$,
we can consider $r$ large enough such that $(1+\epsilon) \xi_k \log r \geq c_r$.
Hence if in round $1 \leq s \leq r$ arm $k$ has sample size of at least $(1+\epsilon) \xi_r \log r$,
it wins against leading optimal arm $\ell$ only if
$$\bar Y_{ku} \geq \bar Y_{\ell,t:(t+u-1)} \mbox{ for some } 1 \leq t \leq n_{\ell}-u+1 (\leq r).
$$
By (\ref{Gkr}), (\ref{4.1}), (\ref{4.2}) and Bonferroni's inequality,
\begin{eqnarray*}
P(G_k^r) & \leq & \sum_{u= \lceil (1+\epsilon) \xi_k \log r \rceil}^{r-1} P \{ \bar Y_{ku} \geq \bar Y_{\ell,t:(t+u-1)} \mbox{ for some } 1 \leq t \leq r \mbox{ and } \ell \in \Xi \} \cr
& \leq & \sum_{u = \lceil (1+\epsilon) \xi_k \log r \rceil}^{r-1} (| \Xi | re^{-u I_1(\omega_k)}+e^{-u I_k(\omega_k)}) \cr
& \leq & \tfrac{Kr}{1-e^{-I_1(\omega_k)}} e^{-(1+\epsilon) \xi_k I_1(\omega_k) \log r} + \tfrac{1}{1-e^{-I_k(\omega_k)}} e^{-(1+\epsilon) \xi_k I_k(\omega_k) \log r},
\end{eqnarray*}
and (I) holds because $(1+\epsilon) \xi_k I_1(\omega_k) >1$ and $(1+\epsilon) \xi_k I_k(\omega_k) >0$.

Let $J_k > \max(\frac{1}{I_k(\omega)},\frac{2}{I_1(\omega)})$.
It follows from (\ref{Hkr}), (\ref{4.1}), (\ref{4.2}) and the arguments above that
\begin{eqnarray*}
P(H_k^r) & \leq & \sum_{u = \lceil J_k \log r \rceil}^{r-1} (| \Xi | re^{-u I_1(\omega)}+e^{-u I_k(\omega)}) \cr
& \leq & \tfrac{Kr}{1-e^{-I_1(\omega)}} e^{-J_k I_1(\omega) \log r}+\tfrac{1}{1-e^{-I_k(\omega)}} e^{-J_k I_k(\omega) \log r},
\end{eqnarray*}
and (II) holds because $J_k I_1(\omega)>2$ and $J_k I_k(\omega)>1$. 
 $\wbox$

\subsection{Proof of Theorem \ref{thm2}}
We consider here SSTC.
By Lemmas \ref{lem2} and \ref{lem1} it suffices,
in Lemmas \ref{lem6}--\ref{lem8} below,
to verify the conditions needed to show that (\ref{xik}) holds with $\xi_k = 1/M(\frac{\mu_*-\mu_k}{\sigma_k})$.
Lemma \ref{lem7a} provides the underlying large deviations bounds for the standard error estimator.
Let $\Phi(z) = P(Z \leq z)$ and $\bar \Phi(z) = P(Z>z) (\leq e^{-z^2/2}$ for $z \geq 0$) for $Z \sim$ N(0,1).

\begin{lem} \label{lem7a}
For $1 \leq k \leq K$ and $t \geq 2$,
\begin{eqnarray} \label{u1}
P(\wht \sigma^2_{kt}/\sigma_k^2 \geq x) & \leq & \exp[\tfrac{(t-1)}{2}(\log x-x+1)] \mbox{ if } x>1, \\ \label{u2}
P(\wht \sigma^2_{kt}/\sigma_k^2 \leq x) & \leq & \exp[\tfrac{(t-1)}{2}(\log x-x+1)] \mbox{ if } 0<x<1. 
\end{eqnarray}
\end{lem}

{\sc Proof}.
We note that $\wht \sigma^2_{kt}/\sigma_k^2 \stackrel{d}{=} \frac{1}{t-1} \sum_{s=1}^{t-1} U_s$,
where $U_s \stackrel{\rm i.i.d.}{\sim} \chi^2_1$,
and that $U_1$ has large deviations rate function
\begin{eqnarray*}
I_U(x) & = & \sup_{\theta<\frac{1}{2}} (\theta x-\log E e^{\theta U_1}) \cr
& = & \sup_{\theta<\frac{1}{2}} [\theta x-\tfrac{1}{2} \log(\tfrac{1}{1-2 \theta})] = \tfrac{1}{2} (x-1-\log x).
\end{eqnarray*}
The last equality holds because the supremum occurs when $\theta=\frac{x-1}{2x}$.
We conclude (\ref{u1}) and (\ref{u2}) from (\ref{LD1}) and (\ref{LD2}) respectively. 
$\wbox$

\begin{lem} \label{lem6}
Under {\rm (\ref{fnormal})}, 
$P(B^r) \leq Q e^{-a r}$ for some $Q>0$ and $a>0$,
when $\frac{r}{K}-1 \geq c_r$. 
\end{lem}

{\sc Proof}.
Let $r$ be such that $\frac{r}{K}-1 \geq c_r$.
The event $B^r$ occurs if at round $r$ the leading arm $\ell$ is optimal, 
and it loses to an inferior arm $k$ with $n_k =u$ and $n_{\ell} = u+1$ for $u \geq \frac{r}{K}-1$.
Let $k \not\in \Xi$,
$\ell \in \Xi$ and let $\epsilon > 0$ be such that $\omega := \frac{\mu_k-\mu_{\ell}+\epsilon}{2 \sigma_k} < 0$.
Let $\tau_i(u)$,
$1 \leq i \leq 3$,
be quantities that we shall define below.
Note that
\begin{equation} \label{4.10}
\tau_1(u) := P(\tfrac{\bar Y_{ku} - \bar Y_{{\ell},u+1}}{\hat \sigma_{ku}} \geq \omega) \leq P(\tfrac{\bar Y_{ku}-\bar Y_{{\ell},u+1}}{2 \sigma_k} \geq \omega) + P(\wht \sigma_{ku} \geq 2 \sigma_k).
\end{equation}
Since $\bar Y_{ku} - \bar Y_{{\ell},u+1} \sim {\rm N}(\mu_k-\mu_{\ell},\tfrac{\sigma_{\ell}^2}{u+1}+\tfrac{\sigma_k^2}{u})$,
\begin{equation} \label{4.8}
P(\tfrac{\bar Y_{ku}-\bar Y_{\ell,u+1}}{2 \sigma_k} \geq \omega) \leq \bar \Phi(\epsilon \sqrt{\tfrac{u}{\sigma_{\ell}^2+\sigma_k^2}})
\leq e^{-\frac{\epsilon^2 u}{2(\sigma_k^2+\sigma_{\ell}^2)}}.
\end{equation}
It follows from (\ref{u1}) and (\ref{u2}) that 
\begin{eqnarray} \label{4.9}
P(\wht \sigma_{ku} \geq 2 \sigma_k) & \leq & e^{-a_1(u-1)/2}, \\ \label{LD2a}
P(\wht \sigma_{\ell u} \leq \tfrac{\sigma_{\ell}}{2}) & \leq & e^{-a_2(u-1)/2}, 
\end{eqnarray}
where $a_1=1-\log 2(>0)$ and $a_2 = \log 2-\frac{1}{2} (>0)$.
By (\ref{4.10})--(\ref{4.9}),
\begin{equation} \label{con1}
\tau_1(u) \leq e^{-\frac{\epsilon^2 u}{2(\sigma_k^2+\sigma_{\ell}^2)}} + e^{-\frac{a_1(u-1)}{2}}.
\end{equation}

Since $\frac{\bar Y_{\ell u}-\bar Y_{\ell,u+1}}{\sigma_{\ell}/2} \sim$ N(0,$\lambda$) for $\lambda \leq 4(\frac{1}{u}+\frac{1}{u+1}) \leq \frac{8}{u}$,
it follows that
$$P(\tfrac{\bar Y_{\ell u}-\bar Y_{\ell,u+1}}{\sigma_{\ell}/2} \leq \omega) \leq \bar \Phi(|\omega| \sqrt{\tfrac{u}{8}}) \leq e^{-\frac{\omega^2 u}{16}}. 
$$
Hence by (\ref{LD2a}),
\begin{eqnarray} \label{4.10a}
\tau_2(u) & := & P(\tfrac{\bar Y_{{\ell},t:(t+u-1)}-\bar Y_{\ell,u+1}}{\hat \sigma_{{\ell},t:(t+u-1)}} \leq \omega \mbox{ for } t=1 \mbox{ or } 2) \\ \nonumber 
& \leq & 2 [P(\tfrac{\bar Y_{{\ell} u}-\bar Y_{\ell,u+1}}{\sigma_{\ell}/2} \leq \omega) + P(\wht \sigma_{\ell u} \leq \tfrac{\sigma_{\ell}}{2})] \\ \nonumber
& \leq & 2(e^{-\frac{\omega^2 u}{16}} + e^{-\frac{a_2(u-1)}{2}}).
\end{eqnarray}
We check that for $\omega_k =\frac{\mu_k+\mu_{\ell}}{2}$,
\begin{eqnarray} \label{klmax}
\tau_3(u) & := & P(\bar Y_{k u} \geq \bar Y_{\ell,u+1}) \\ \nonumber
& \leq & P(\bar Y_{k u} \geq \omega_k) + P(\bar Y_{\ell,u+1} \leq \omega_k) \\ \nonumber
& \leq & e^{-u(\omega_k-\mu_k)^2/(2 \sigma_k^2)} + e^{-(u+1)(\omega_k-\mu_{\ell})^2/(2 \sigma_{\ell}^2)}.
\end{eqnarray}

By (\ref{con1})--(\ref{klmax}),
$$P(B^r) \leq  \sum_{k \not\in \Xi} \sum_{\ell \in \Xi} \sum_{u = \lceil \frac{r}{K} \rceil-1}^r [\tau_1(u)+\tau_2(u)+\tau_3(u)],  
$$
and Lemma \ref{lem6} indeed holds. 
$\wbox$

\begin{lem} \label{lem7}
Under {\rm (\ref{fnormal})}, 
$P(C^r) \leq K^2 e^{-c_r a} \tfrac{(\log r)^6}{r} + o(r^{-1})$ for some $a>0$.  
\end{lem}

{\sc Proof}.
The event $C^r$ occurs if at round $r$ the leading arm $k$ is inferior,
and it wins a challenge against one or more optimal arms $\ell$.
By step 2(b)ii. of SSTC,
we need only consider  $n_{\ell} \geq c_n$.
Note that when $n_k =n_{\ell}$,
for arm $k$ to be leader,
by the tie-breaking rule we require $\bar Y_{k n_k} \geq \bar Y_{\ell n_{\ell}}$.
Consider $n_k$ taking values $u$, 
$n_{\ell}$ taking values $v$ and let $\tau_i(\cdot)$,
$1 \leq i \leq 4$, 
be quantities that we shall define below.

Case 1. $n_{\ell} > (\log r)^2$.
Let $\omega = \frac{\mu_{\ell}+\mu_k}{2}$ and check that
\begin{eqnarray} \label{7.1}
\tau_1(u,v) & := & P(\bar Y_{\ell v} \leq \omega) + P(\bar Y_{ku} \geq \omega) \\ \nonumber
& \leq & e^{-v(\mu_{\ell}-\mu_k)^2/(8 \sigma_{\ell}^2)} + e^{-u(\mu_{\ell}-\mu_k)^2/(8 \sigma_k^2)}. 
\end{eqnarray}

Case 2. $(c_r \leq) c_n \leq n_{\ell} < (\log r)^2$. 
Let $\omega$ be such that 
\begin{equation} \label{Pomega}
(p_{\omega}:=)P( \tfrac{\bar Y_{kv} - \mu_k+r^{-\frac{1}{3}}}{\hat \sigma_{kv}} \leq \omega) = \tfrac{(\log r)^4}{r}.
\end{equation}
Hence
\begin{eqnarray} \label{pY}
\tau_2(v) & := & P( \tfrac{\bar Y_{k,t:(t+v-1)} - \mu_k+r^{-\frac{1}{3}}}{\hat \sigma_{kv}} > \omega \mbox{ for } 1 \leq t \leq \tfrac{r}{K}) \\ \nonumber
& [\leq & (1-p_{\omega})^{\lfloor \frac{r}{Kv} \rfloor}] \leq \exp[-\tfrac{(\log r)^2}{K}+1].
\end{eqnarray}

We shall show that there exists $a>0$ such that for large $r$,
\begin{equation} \label{7.32}
\tau_3(v) := P(\tfrac{\bar Y_{\ell v}-\mu_k-r^{-\frac{1}{3}}}{\hat \sigma_{\ell v}} \leq \omega) \leq \tfrac{e^{-av} (\log r)^4}{r}
(\leq \tfrac{e^{-c_r a} (\log r)^4}{r}).
\end{equation}
For $u \geq \tfrac{r}{K}$,
\begin{equation} \label{tau4}
\tau_4(u) := P(|\bar Y_{ku}-\mu_k| \geq r^{-\frac{1}{3}}) \leq e^{-u r^{-1/3}/(2 \sigma_k^2)} \leq e^{-r^{2/3}/(2K \sigma_k^2)}.
\end{equation}
Since (\ref{pY}) and (\ref{7.32}) hold with ``$-\bar Y_{ku}$" replacing ``$-\mu_k+r^{-\frac{1}{3}}$" and ``$-\mu_k-r^{-\frac{1}{3}}$" respectively,
by adding $\tau_4(u)$ to the upper bounds,
$$P(C^r) \leq \sum_{k \not\in \Xi} \sum_{\ell \in \Xi} \Big( \sum_{v = \lceil c_r \rceil}^{\lfloor (\log r)^2 \rfloor} [\tau_2(v)+\tau_3(v)] + 
\sum_{u=\lceil \frac{r}{K} \rceil}^r 2 \tau_4(u)  + \sum_{u = \lceil \frac{r}{K} \rceil}^{r} \sum_{v = \lceil (\log r)^2 \rceil}^r \tau_1(u,v) \Big).
$$
We conclude Lemma \ref{lem7} from (\ref{7.1}) and (\ref{pY})--(\ref{tau4}).

We shall now show (\ref{7.32}),
noting firstly that for $r$ large,
the $\omega$ satisfying (\ref{Pomega}) is negative.
This is because for $v < (\log r)^2$,
\begin{eqnarray*}
P(\tfrac{\bar Y_{kv}-\mu_k+r^{-\frac{1}{3}}}{\hat \sigma_{kv}} \leq 0) = \Phi(-\tfrac{r^{-\frac{1}{3}} \sqrt{v}}{\sigma_k}) \rightarrow \tfrac{1}{2},
\end{eqnarray*}
whereas $\frac{(\log r)^4}{r} \rightarrow 0$. 

Let $g_v$ be the common density function of $\wht \sigma_{kv}/\sigma_k$ and $\wht \sigma_{\ell v}/\sigma_{\ell}$.
By the independence of $\bar Y_{k v}$ and $\wht \sigma_{k v}$,
\begin{eqnarray} \label{g1}
P(\tfrac{\bar Y_{k v}-\mu_k+r^{-\frac{1}{3}}}{\hat \sigma_{k v}} \leq \omega) & = & 
\int_0^{\infty} P(\tfrac{\bar Y_{k v}-\mu_k+r^{-\frac{1}{3}}}{\sigma_k} \leq \omega x) g_v(x) dx \\ \nonumber
& = & \int_0^{\infty} \Phi(\sqrt{v}(\omega x-\tfrac{r^{-\frac{1}{3}}}{\sigma_k})) g_v(x) dx.
\end{eqnarray}

By similar arguments, 
\begin{equation} \label{g2}
P(\tfrac{\bar Y_{\ell v}-\mu_k-r^{-\frac{1}{3}}}{\hat \sigma_{\ell v}} \leq \omega) 
= \int_0^{\infty} \Phi(\sqrt{v}(\omega x-\tfrac{\Delta-r^{-\frac{1}{3}}}{\sigma_{\ell}})) g_v(x) dx,
\end{equation}
where $\Delta := \mu_{\ell}-\mu_k(>0)$.
Let $\delta_1 = \frac{r^{-\frac{1}{3}}}{\sigma_k}$,
$\delta_2 = \frac{\Delta - r^{-\frac{1}{3}}}{\sigma_{\ell}}$ and $b=-\omega x$.
Since $b>0$ and $\delta_2 > \delta_1 > 0$ for $r$ large, 
\begin{equation} \label{Phin}
\Phi(\sqrt{v}(-b-\delta_2)) \leq e^{-a_r v}\Phi(\sqrt{v}(-b-\delta_1)),
\end{equation}
where $a_r = \frac{(\delta_2-\delta_1)^2}{2}$ ($\rightarrow \frac{\Delta^2}{2 \sigma_{\ell}^2}$ as $r \rightarrow \infty$).
Let $a=\frac{\Delta_2^2}{4 \sigma_{\ell}^2}$. 
It follows from (\ref{g1})--(\ref{Phin}) that for $r$ large,
$$P(\tfrac{\bar Y_{\ell v}-\mu_k-r^{-\frac{1}{3}}}{\hat \sigma_{kv}} \leq \omega) \leq e^{-av} P(\tfrac{\bar Y_{kv}-\mu_k+r^{-\frac{1}{3}}}{\hat \sigma_{kv}} \leq \omega).
$$
Hence by (\ref{Pomega}), 
the inequality in (\ref{7.32}) indeed holds.  
$\wbox$

\begin{lem} \label{lem10a}
Let $Z_s \sim {\rm N}(0,\frac{1}{s+1})$ and $W_s \sim \chi^2_s/s$ be independent.
For any $g<0$ and $0 < \delta < M(g)$,
there exists $Q>0$ such that for $s_1 \geq 1$,
$$\sum_{s=s_1}^{\infty} P \{ \tfrac{Z_s}{\sqrt{W_s}} \leq g \} \leq Q e^{-s_1[M(g)-\delta]}.
$$
\end{lem}

{\sc Proof}.
Consider the domain $\Omega = {\bf R}^+ \times {\bf R}$, 
and the set
$$A = \{ (w,z) \in \Omega: z \leq g \sqrt{\omega} \}.
$$
Let $I(w,z) = \frac{1}{2}(z^2+w-1-\log w)$,
and check that
\begin{eqnarray} \label{Iwz}
& & \inf_{(w,z) \in A} I(w,z) = \inf_{w>0} I(w,g \sqrt{w}) \\ \nonumber
& = & \inf_{w>0} [\tfrac{1}{2}(g^2 w+w-1-\log w)] = \tfrac{1}{2} \log(1+g^2) = M(g),
\end{eqnarray}
the second last equality follows from the infimum occurring at $w=\frac{1}{g^2+1}$.

Let $L_v$,
$1 \leq v \leq V$,
be half-spaces constructed as follows. 
Let 
\begin{eqnarray} \label{H1}
L_1 & = & \{ (w,z): z \leq z_1, 0 < w < \infty \}, \\ \nonumber
& & \mbox{ with } g < z_1 < 0 \mbox{ and } I(1,z_1) \geq M(g)-\delta. 
\end{eqnarray}
The existence of $z_1$ satisfying second line of (\ref{H1}) follows from $I(1,g) = \frac{1}{2}g^2 > M(g)$.
Since $(A \setminus L_1) \subset (0,1) \times (z_1,0)$,
by (\ref{Iwz}), 
we can find half-spaces
\begin{eqnarray} \label{Hm}
L_v & = & \{ (w,z): 0 < w \leq w_v, z \leq z_v \} \mbox{ with } 0  < w_v < 1, \\ \nonumber
& & \quad z_v \leq 0 \mbox{ and } I(w_v,z_v) \geq M(g)-\delta, \quad 2 \leq v \leq V,
\end{eqnarray} 
such that $(A \setminus L_1) \subset \cup_{v=2}^V L_v$.
Therefore $A \subset \cup_{v=1}^V L_v$,
and so 
\begin{equation} \label{s0in}
\sum_{s=s_1}^{\infty} P \{ \tfrac{Z_s}{\sqrt{W_s}} \leq g \} \leq \sum_{s=s_1}^{\infty} \sum_{v=1}^V P \{ (W_s,Z_s) \in L_v \}. 
\end{equation}

It follows from (\ref{u2}), (\ref{H1}), (\ref{Hm}) and the independence of $Z_s$ and $W_s$,
setting $w_1=1$, 
that
\begin{equation} \label{WsZs}
P \{ (W_s,Z_s) \in L_v \} \leq e^{-s I(w_v,z_v)} \leq e^{-s[M(g)-\delta]}, \ 1 \leq v \leq V.
\end{equation}
Lemma \ref{lem10a},
with $Q=\frac{V}{1-e^{-M(g)+\delta}}$,
follows from substituting (\ref{WsZs}) into (\ref{s0in}). 
$\wbox$

\begin{lem} \label{lem8}
Under {\rm (\ref{fnormal})} and $c_r=o(\log r)$,
{\rm (I)} $($in the statement of Lemma~${\ref{lem2}})$ holds for $\xi_k = 1/M(\tfrac{\mu_*-\mu_k}{\sigma_k})$ and {\rm (II)} holds for $J_k$ large.
\end{lem}

{\sc Proof}. 
By considering the rewards $Y_{kt}-\mu_*$,
we may assume without loss of generality that $\mu_*=0$.
Let $k \not\in \Xi$ (hence $\mu_k<0$) and $\epsilon>0$.
Let $g_k = \frac{\mu_k}{\sigma_k}$ and let $g_{\omega} < 0$ and $\delta>0$ be such that
\begin{equation} \label{gomega}
0 > g_{\omega}-3 \delta > g_k \mbox{ and } (1+\epsilon) [M(g_{\omega}-\delta)-\delta] > M(g_k).
\end{equation}
Let $m_r = \lceil (1+\epsilon) (\log r)/M(g_k) \rceil$.
Since $c_r=o(\log r)$,
we can consider $r$ large enough such that $m_r \geq c_r$.
By (\ref{u2}), 
\begin{equation} \label{sdbound}
\sum_{u = m_r}^r P(\wht \sigma^2_{\ell u}/\sigma_{\ell}^2 \leq \tfrac{1}{4}) \rightarrow 0, \quad 1 \leq \ell \leq K.
\end{equation}
Let $\sigma_0 = \min_{1 \leq \ell \leq K} \sigma_{\ell}$.
For $\ell \in \Xi$,
\begin{eqnarray} \label{7b}
& & \sum_{v = \lceil \frac{r}{K} \rceil}^r P(\tfrac{|\bar Y_{\ell v}|}{\sigma_0/2} \geq \delta)  
\leq \sum_{v = \lceil \frac{r}{K} \rceil}^r \exp (-\tfrac{\delta^2 \sigma_0^2 v}{8 \sigma_{\ell}^2} \Big) =O(r^{-1}), \\ \label{YY}
& & \eta^r := P(\bar Y_{k n_k} \geq \bar Y_{\ell n_{\ell}} \mbox{ for some } n_k \geq m_r, n_{\ell} \geq \tfrac{r}{K}, \ell \in \Xi) \\ \nonumber
& & \qquad \leq \sum_{u=m_r}^r \exp(-\tfrac{u \mu_k^2}{8 \sigma_k^2})
+ \sum_{\ell \in \Xi} \sum_{v = \lceil \frac{r}{K} \rceil}^r \exp(-\tfrac{v \mu_k^2}{8 \sigma_{\ell}^2}) \rightarrow 0.
\end{eqnarray}

By (\ref{Gkr}) and (\ref{gomega}),
\begin{eqnarray} \label{PGkr}
& & P(G_k^r) \leq P(\tfrac{\bar Y_{kn_k}-\bar Y_{\ell n_{\ell}}}{\hat \sigma_{k n_k}} \geq \tfrac{\bar Y_{\ell,t:(t+n_k-1)}-\bar Y_{\ell n_{\ell}}}{\hat \sigma_{\ell n_k}} \\ \nonumber
& & \quad \mbox{ for some } 1 \leq t \leq r, \ell \in \Xi, n_k \geq m_r, n_{\ell} \geq \tfrac{r}{K})+\eta^r \\ \nonumber
& & \qquad \quad \leq \sum_{u=m_r}^r \Big[ P(\tfrac{\bar Y_{ku}}{\hat \sigma_{ku}} \geq g_k+\delta) 
+ r \sum_{\ell \in \Xi} P(\tfrac{\bar Y_{\ell u}}{\hat \sigma_{\ell u}} \leq g_{\omega}-\delta)
\\ \nonumber
& & \qquad \qquad \quad +\sum_{\ell=1}^K P(\wht \sigma^2_{\ell u}/\sigma^2_{\ell} \leq \tfrac{1}{4}) \Big] 
+ \sum_{\ell \in \Xi} \sum_{v = \lceil \frac{r}{K} \rceil}^r P(\tfrac{|\bar Y_{\ell v}|}{\sigma_0/2} \geq \delta) + \eta^r.
\end{eqnarray}

By (\ref{sdbound})--(\ref{PGkr}),
to show (I),
it suffices to show that
\begin{eqnarray} \label{show1}
\sum_{u=m_r}^r P(\tfrac{\bar Y_{k u}}{\hat \sigma_{k u}} \geq g_k+\delta) & \rightarrow & 0, \\ \label{show3}
r \sum_{u=m_r}^r P(\tfrac{\bar Y_{\ell u}}{\hat \sigma_{\ell u}} \leq g_{\omega}-\delta) & \rightarrow & 0.
\end{eqnarray}

Keeping in mind that $g_k+\delta < 0$, 
let $w>1$ be such that $\sqrt{w}(g_k+\delta) > g_k$.
It follows from (\ref{u1}) and $g_k \sigma_k = \mu_k$ that
\begin{eqnarray*}
& & \sum_{u=m_r}^r P(\tfrac{\bar Y_{ku}}{\hat \sigma_{k u}} \geq g_k+\delta) \cr
& \leq & \sum_{u=m_r}^r [P(\bar Y_{ku} \geq \sqrt{w}(\mu_k+\delta \sigma_k))+P(\wht \sigma^2_{ku}/\sigma_k^2 \geq w)] \cr
& \leq & \sum_{u=m_r}^r [e^{-u[\mu_k-\sqrt{w}(\mu_k+\delta \sigma_k)]^2/(2 \sigma_k^2)}+e^{-(u-1)(w+1-\log w)/2}], 
\end{eqnarray*}
and (\ref{show1}) indeed holds.
Finally by Lemma \ref{lem10a},
$$\sum_{u=m_r}^r P(\tfrac{\bar Y_{\ell u}}{\hat \sigma_{\ell u}} \leq g_{\omega}-\delta) \leq Q e^{-(m_r-1)[M(g_{\omega}-\delta)-\delta]},
$$
for some $Q>0$,
and so (\ref{show3}) follows from (\ref{gomega}). 

To show (II),
we consider $m_r = \lceil J_r \log r \rceil$.
By (\ref{u2}),
we can select $J_k$ large enough to satisfy (\ref{sdbound}) with ``$\rightarrow 0$" replaced by ``$=O(r^{-1})$". 
We note that (\ref{PGkr}) holds with $H_k^r$ in place of $G_k^r$ for this $m_r$.
Therefore to show (II),
it suffices to note that for $J_k$ large enough,
(\ref{YY}), 
(\ref{show1}) and (\ref{show3}) hold with ``$\rightarrow 0$" replaced by ``$=O(r^{-1})$".
$\wbox$

\subsection{Proof of Theorem \ref{thm3}}

Assume (C1)--(C3) and let $\wtd \mu = \max_{k: \mu_k < \mu^*} \mu_k$.
By Lemmas \ref{lem2} and \ref{lem1} it suffices,
in Lemmas \ref{lem9}--\ref{lem11} below,
to verify the conditions needed for SSMC to satisfy (\ref{xik}) for some $\xi_k>0$.

\begin{lem} \label{lem9}
Under {\rm (C2)},
$P(B^r) \leq \frac{3QK^2}{1-e^{-b}} e^{-b(\frac{r}{K}-1)}$ for some $b>0$ and $Q>0$,
when $\frac{r}{K}-1 \geq c_r$.
\end{lem}

{\sc Proof}. 
Consider $r$ such that $(n_k \geq) \frac{r}{K}-1 \geq c_r$.
Let $\epsilon = \frac{1}{2}(\mu_*-\wtd \mu)$ and let $b$ and $Q$ be the constants satisfying (C2). 
Lemma \ref{lem9} follows from arguments similar to those in the proof of Lemma \ref{lem3a},
setting $\omega = \frac{1}{2} (\mu_*+\wtd \mu)$.
$\wbox$

\begin{lem} \label{lem10}
Under {\rm (C1)--(C3)},
$P(C^r) \leq K^2 Q_1 e^{-c_r b_1} \frac{(\log r)^6}{r} +o(r^{-1})$ for some $b_1>0$ and $Q_1>0$.
\end{lem}

{\sc Proof}.
The event $C^r$ occurs if at round $r$ the leading arm $k$ is inferior, 
and it wins against one or more optimal arms $\ell$.
By step 2(b)ii. of SSMC,
we need only consider $n_{\ell} =v$ for $v \geq c_n$.
Note that $n_k \geq \frac{r}{K}$ and $n_k \geq n_{\ell}$.

Case 1: $n_{\ell} > (\log r)^2$.
Let $\omega$ and $\epsilon$ be as in the proof of Lemma \ref{lem9}.
By (C2),
there exists $b>0$ and $Q>0$ such that
\begin{equation} \label{case1}
P(\bar Y_{\ell v} \leq \omega)+P(\bar Y_{kv} \geq \omega) \leq 2Qe^{-vb}.
\end{equation}

Case 2: $n_{\ell}=v$ for $(c_r \leq) c_n \leq v < (\log r)^2$.
Select $\omega(\leq \mu_k$ for $r$ large) such that
\begin{equation} \label{bark}
P(\bar Y_{kv} < \omega) \leq \tfrac{(\log r)^4}{r} \leq P(\bar Y_{kv} \leq \omega).
\end{equation}
Let $p_{\omega} = P(\bar Y_{kv} > \omega)$ and let $d= \lceil 2 (\log r)^2 \rceil$, 
$\eta = \lfloor \frac{r/K-1}{d} \rfloor$.
By (C1) and the second inequality of (\ref{bark}),
\begin{eqnarray} \label{pstar}
\tau(v)& := & P(\bar Y_{k,t:(t+v-1)} > \omega \mbox{ for } 1 \leq t \leq \tfrac{r}{K}) \\ \nonumber
& \leq & P(\bar Y_{k,t:(t+v-1)} > \omega \mbox{ for } t=1,d+1,\ldots, \eta d+1) \\ \nonumber
& \leq & p_{\omega}^{\eta+1} + \eta [1-\lambda_k({\bf R})]^{d-v+1} \\ \nonumber
& \leq & \exp(-\tfrac{(\eta+1)(\log r)^4}{r}) + \eta [1-\lambda_k({\bf R})]^{(\log r)^2} [=o(r^{-2})].
\end{eqnarray}

To see the second inequality of (\ref{pstar}),
let
$$D_m = \{ \bar Y_{k,t:(t+v-1)} > \omega \mbox{ for } t=md+1 \}, \quad 0 \leq m \leq \eta.
$$
Note that the probability in the second line of (\ref{pstar}) is $P(\cap_{m=0}^{\eta} D_m)$,
and that by (\ref{bark}),
$P(D_m) =p_{\omega} \leq 1-\frac{(\log r)^4}{r}$.
By the triangular inequality and the convention $\prod_{m=\eta+1}^{\eta}=1$,
\begin{eqnarray} \label{tri}
& & \Big| P(\cap_{m=0}^{\eta} D_m) - \prod_{m=0}^{\eta} P(D_m) \Big| \\ \nonumber
& \leq & \sum_{u=1}^{\eta} \Big| P(\cap_{m=0}^u D_m) \prod_{m=u+1}^{\eta} P(D_m) - P(\cap_{m=0}^{u-1} D_m) \prod_{m=u}^{\eta} P(D_m) \Big| \\ \nonumber
& \leq & \sum_{u=1}^{\eta} |P(\cap_{m=0}^u D_m) - P(\cap_{m=0}^{u-1} D_m) P(D_u)|.
\end{eqnarray}
By (C1),
\begin{equation} \label{Dm}
|P(\cap_{m=0}^u D_m) - P(\cap_{m=0}^{u-1} D_m) P(D_u)| \leq [1-\lambda_k({\bf R})]^{d-v+1}, \quad 1 \leq u \leq \eta,
\end{equation}
since $\cap_{m=0}^{u-1} D_m$ depends on $(Y_{k1}, \ldots, Y_{k,(u-1)d+v})$ whereas $D_u$ depends on $(Y_{k,ud+1}, \ldots, Y_{k,ud+v})$.
Substituting (\ref{Dm}) into (\ref{tri}) gives us the second inequality of (\ref{pstar}).

It follows from (C3) and the first inequality of (\ref{bark}) that there exists $Q_1>0$, 
$b_1>0$ and $t_1 \geq 1$ such that for $v \geq t_1$,
$$P(\bar Y_{\ell v} < \omega) \leq Q_1 e^{-b_1 v} \tfrac{(\log r)^4}{r}.
$$
Hence by (\ref{case1}) and (\ref{pstar}),
for $r$ such that $c_r \geq t_1$,
$$P(C^r) \leq \sum_{k \not\in \Xi} \sum_{\ell \in \Xi} \Big( \sum_{v = \lceil (\log r)^2 \rceil}^r 2Q e^{-v b} 
+ \sum_{v = \lceil c_r \rceil}^{\lfloor (\log r)^2 \rfloor} [Q_1 e^{-b_1 c_r} \tfrac{(\log r)^4}{r} + \tau(v)] \Big),
$$
and Lemma \ref{lem10} holds.
$\wbox$

\begin{lem} \label{lem11}
Under {\rm (C2)} and $c_r=o(\log r)$,
statement {\rm (II)} in Lemma $\ref{lem2}$ holds.
\end{lem}

{\sc Proof}. 
Let $\epsilon$ and $\omega$ be as in the proof of Lemma \ref{lem9},
and let $b$ and $Q$ be the constants satisfying (C2).
For an optimal arm $\ell$,
\begin{eqnarray*}
P(\bar Y_{{\ell},t:(t+u-1)} \leq \omega \mbox{ for some } 1 \leq t \leq r) & \leq & Qre^{-ub}, \cr
P(\bar Y_{ku} \geq \omega) & \leq & Qe^{-ub}.
\end{eqnarray*}
Let $J_k > \frac{2}{b}$.
Since $c_r=o(\log r)$,
for $r$ large, 
$\lceil J_k \log r \rceil \geq c_r$ and therefore by Bonferroni's inequality,
$$P(H_k^r) \leq \sum_{\ell \in \Xi} \sum_{u = \lceil J_k \log r \rceil}^r Q(r+1) e^{-ub},
$$ 
and (II) holds.
$\wbox$

\begin{appendix}
\section{Showing (\ref{min1})}

Let $\Phi(z) = P(Z \leq z)$ for $Z \sim N(0,1)$. 
It follows from $\Phi(-z) = [1+o(1)] \frac{1}{z \sqrt{2 \pi}} e^{-z^2/2}$ as $z \rightarrow \infty$ that
\begin{eqnarray} \label{Mill1}
\Phi(-\sqrt{2 \log n}) & = & \tfrac{1+o(1)}{2n \sqrt{\pi \log n}}, \\ \label{Mill2}
\Phi(-\sqrt{2 \log (\tfrac{n}{(\log n)^2})}) & = & [1+o(1)] \tfrac{(\log n)^{3/2}}{2n \sqrt{\pi}}. 
\end{eqnarray}
Assume without loss of generality $\mu_1=0$ and consider $n_1=u$ and $n_2=v$ (hence $u+v=n$) with $v = O(\log n)$.
By (\ref{Mill1}) and Bonferroni's inequality,
\begin{eqnarray} \label{bound1}
& & P(\min_{1 \leq t \leq u-v+1} \bar Y_{1,t:(t+v-1)} \leq - \sqrt{\tfrac{2 \log n}{v}}) \\ \nonumber
& \leq & \sum_{t=1}^{u-v+1} P(\bar Y_{1,t:(t+v-1)} \leq -\sqrt{\tfrac{2 \log n}{v}}) \\ \nonumber
& = & (u-v+1) \Phi(-\sqrt{2 \log n}) \rightarrow 0.
\end{eqnarray}
By (\ref{Mill2}) and independence of $\bar Y_{1,(sv+1):[(s+1)v]}$ for $0 \leq s \leq \frac{u-v}{v}$,
\begin{eqnarray} \label{bound2}
& & P(\min_{1 \leq t \leq u-v+1} \bar Y_{1,t:(t+v-1)} \geq - \sqrt{\tfrac{2 \log(n/(\log n)^2)}{v}}) \\ \nonumber
& \leq & P(\min_{0 \leq s \leq (u-v)/v} \bar Y_{1,(sv+1):[(s+1)v]} \geq -\sqrt{\tfrac{2 \log(n/(\log n)^2)}{v}}) \\ \nonumber
& = & [1-\Phi(-\sqrt{2 \log (\tfrac{n}{(\log n)^2})})]^{\lfloor \frac{u-v}{v} \rfloor+1} \\ \nonumber
& \leq & \exp[-(\lfloor \tfrac{u-v}{v} \rfloor+1) \Phi(-\sqrt{2 \log (n/(\log n)^2)})] \rightarrow 0.
\end{eqnarray}
We conclude (\ref{min1}) from (\ref{bound1}) and (\ref{bound2}).

\section{Verifications of (C1)--(C3) for double exponential densities}

By dividing $Y_{kt}$ by $\tau$ if necessary,
we may assume without loss of generality that $\tau=1$.
We check that (C1) holds for $\lambda_k(A) = \int_A f_k(y) dy$,
whereas (C2) follows from the Chernoff bounds given in Lemma \ref{lem3b},
that is (\ref{YQ}) holds for $Q=2$ and $b = I(\epsilon)$,
where $I(\mu) = \sup_{|\theta|<1} [\theta \mu-\log(1-\theta^2)]$ is the large deviations rate function of the double exponential density $f(y)=\frac{1}{2} e^{-|y|}$.

Let $S_t = \sum_{u=1}^t Y_u$ with $Y_u \stackrel{\rm i.i.d.}{\sim} f$ and let $\Delta=\mu_{\ell}-\mu_k$.
Since $\mu_k-Y_{kt} \sim f$,
and similarly when $k$ is replaced by $\ell$,
to show (C3),
it suffices to show that for $z \geq 0$ and $t \geq 1$,
\begin{equation} \label{B2}
P(S_t > z+\Delta t) \leq e^{-tb_1} P(S_t > z),
\end{equation}
where $b_1=\Delta-2 \log(1+\frac{\Delta}{2}) (>0)$.
By (\ref{B2}),
(C3) holds for $Q_1=1$,
$t_1=1$ and the above $b_1$.

Since $Y_u \stackrel{d}{=} Z_{u1}-Z_{u2}$,
with $Z_{u1}$ and $Z_{u2}$ independent exponential random variables with mean 1,
it follows that $S_t \stackrel{d}{=} S_{t1}-S_{t2}$ where $S_{t1}$ and $S_{t2}$ are independent Gamma random variables.
Using this,
Kotz, Kozubowski and Pod\'{g}orski (2001) showed,
see their (2.3.25),
that the density $f_t$ of $S_t$ can be expressed as $f_t(x) = e^{-x} g_t(x)$ for $x \geq 0$, 
where
\begin{equation} \label{coeff}
g_t(x) = \tfrac{1}{(t-1)! 2^{2t-1}} \sum_{j=0}^{t-1} c_{tj} x^j, \mbox{ with } c_{tj} = \tfrac{(2t-2-j)! 2^j}{j! (t-1-j)!}.
\end{equation}

We shall show that
\begin{equation} \label{f2}
g_t'(x)(1+\tfrac{x}{2t}) \leq g_t(x).
\end{equation}
By (\ref{f2}),
$$\tfrac{f_t'(x)}{f_t(x)} = \tfrac{e^{-x}[g_t'(x)-g_t(x)]}{e^{-x} g_t(x)} \leq \tfrac{2t}{x+2t}-1,
$$
and therefore for $y \geq 0$,
\begin{eqnarray*}
\log [\tfrac{f_t(y+t \Delta)}{f_t(y)}] & = & \int_y^{y+ t\Delta} \tfrac{f_t'(x)}{f_t(x)} dx \cr
& \leq & 2t \log(\tfrac{y+(2+\Delta)t}{y+2t}) - t \Delta \leq -tb_1.
\end{eqnarray*}
Hence $f_t(y+t \Delta) \leq e^{-tb_1} f_t(y)$.
It follows that for $z \geq 0$,
\begin{eqnarray*}
P(S_t > z+t \Delta) & = & \int_z^{\infty} f_t(y+t \Delta) dy \cr
& \leq & e^{-tb_1} \int_z^{\infty} f_t(y) dy  = e^{-tb_1} P(S_t>z), 
\end{eqnarray*}
and (C3) indeed holds.

We shall now show (\ref{f2}) by checking that after substituting (\ref{coeff}) into (\ref{f2}),
the coefficient of $x^j$ in the left-hand side of (\ref{f2}) is not more than in the right-hand side,
for $0 \leq j \leq t-1$.
More specifically that (with $c_{tt}=0$),
\begin{equation} \label{ctj}
(j+1) c_{t,j+1} + \tfrac{j}{2t} c_{tj} \leq c_{tj} [\Leftrightarrow c_{t,j+1} \leq \tfrac{1}{j+1}(1-\tfrac{j}{2t}) c_{tj}].
\end{equation}
Indeed by (\ref{coeff}),
$$c_{t,j+1} = \tfrac{2(t-1-j)}{(j+1)(2t-2-j)} c_{tj} = \tfrac{1}{j+1}(1-\tfrac{j}{2t-2-j}) c_{tj}, \ 0 \leq j \leq t-1,
$$
and the right-inequality of (\ref{ctj}) holds.
\end{appendix}

\section*{Acknowledgment}
We would like to thank three referees and an Associate Editor for going over the manuscript carefully, 
and providing useful feedbacks.
The changes made in response to their comments have resulted in a much better paper.
Thanks also to Shouri Hu for going over the proofs and performing some of the simulations in Examples 5 and 6.


\begin{thebibliography}{9}
\bibitem{Agr95}
\textsc{Agrawal, R.} (1995).
Sample mean based index policies with $O(\log n)$ regret for the multi-armed bandit problem.
\textit{Adv. Appl. Probab.} \textbf{17} 1054--1078.

\bibitem{ATA89}
\textsc{Agrawal, R., Teneketzis, D.} and \textsc{Anantharam, V.} (1989).
Asymptotically efficient adaptive allocation schemes for controlled Markov chains: Finite parameter space.
\textit{IEEE Trans. Automat. Control} \textbf{AC-34} 1249--1259.

\bibitem{ACF02}
\textsc{Auer, P., Cesa-Bianchi, N.} and \textsc{Fischer, P.} (2002).
Finite-time analysis of the multiarmed bandit problem.
\textit{Machine Learning} \textbf{47} 235--256.

\bibitem{BMM14}
\textsc{Baransi, A., Maillard, O.A.} and \textsc{Mannor, S.} (2014).
Sub-sampling for multi-armed bandits.
\textit{Proceedings of the European Conference on Machine Learning} pp.13.

\bibitem{BF85}
\textsc{Berry, D.} and \textsc{Fristedt, B.} (1985).
\textit{Bandit problems}. Chapman and Hall, London. 

\bibitem{BL02}
\textsc{Brezzi, M.} and \textsc{Lai, T.L.} (2002).
Optimal learning and experimentation in bandit problems.
\textit{J. Econ. Dynamics Cont.} \textbf{27} 87--108.

\bibitem{BK96}
\textsc{Burnetas, A.} and \textsc{Katehakis, M.} (1996).
Optimal adaptive policies for sequential allocation problems.
\textit{Adv. Appl. Math.} \textbf{17} 122--142.

\bibitem{BLL15}
\textsc{Burtini, G., Loeppky, J.} and \textsc{Lawrence, R.} (2015).
A survey of online experiment design with the stochastic multi-armed bandit.
arXiv:1510.00757.

\bibitem{CGMMS13}
\textsc{Capp\'{e}, O., Garivier, A., Maillard, J., Munos, R., Stoltz, G.} (2013).
Kullback-Leibler upper confidence bounds for optimal sequential allocation.
\textit{Ann. Statist.} \textbf{41} 1516--1541.

\bibitem{CL87}
\textsc{Chang, F.} and \textsc{Lai, T.L.} (1987).
Optimal stopping and dynamic allocation.
\textit{Adv. Appl. Probab.} \textbf{19} 829--853.

\bibitem{Git79}
\textsc{Gittins, J.C.} (1979). 
Bandit processes and dynamic allocation indices.
\textit{JRSS`B'} \textbf{41} 148--177.

\bibitem{GJ79}
\textsc{Gittins, J.C.} and \textsc{Jones, D.M.} (1979).
A dynamic allocation index for the discounted multi-armed bandit problem.
\textit{Biometrika} \textbf{66} 561--565.

\bibitem{GL97}
\textsc{Graves, L.} and \textsc{Lai, T.L.} (1997).
Asymptotically efficient adaptive choices of control laws in controlled Markov chains.
\textit{SIAM Control Optim.} \textsc{35} 715--743.


%\bibitem{HW89}
%\textsc{Hu, I.} and \textsc{Wei, C.Z.} (1989).
%Irreversible adaptive allocation rules.
%\textit{Ann. Statist.} \textbf{17} 801--823.

\bibitem{Kau14}
\textsc{Kaufmann, E.} (2014).
Analyse de strat\'{e}gies bay\'{e}siennes et fr\'{e}quentistes pour l'allocation s\'{e}quentielle de ressources,
PhD thesis.

\bibitem{KCG12}
\textsc{Kaufmann, E., Capp\'{e}} and \textsc{Garivier, A.} (2012).
On Bayesian upper confidence bounds for bandit problems.
\textit{Proceedings of the Fifteenth International Conference on Artificial Intelligence and Statistics}
\textbf{22} 592--600.

%\bibitem{KL16}
%\textsc{Kim, M.} and \textsc{Lim, A.} (2016).
%Robust multiarmed bandit problems.
%\textit{Management Sci.} \textbf{62} 264--285.

\bibitem{KKM13}
\textsc{Korda, N., Kaufmann, E.} and \textsc{Munos, R.} (2013).
Thompson sampling for 1-dimensional exponential family bandits.
\textit{NIPS} \textbf{26} 1448--1456.

\bibitem{KKP01}
\textsc{Kotz, S., Kozubowski, T.} and \textsc{Podg\'{o}rski, K.} (2001).
\textit{The Laplace Distribution with Generalizations}, Springer.

\bibitem{KP00}
\textsc{Kuleshov, V.} and \textsc{Precup, D.} (2014).
Algorithms for the multi-armed bandit problem.
arXiv:1402.6028.

\bibitem{Lai87}
\textsc{Lai, T.L.} (1987).
Adaptive treatment allocation and the multi-armed bandit problem.
\textit{Ann. Statist.} \textbf{15} 1091--1114.

\bibitem{LR85}
\textsc{Lai, T.L.} and \textsc{Robbins, H.} (1985).
Asymptotically efficient adaptive allocation rules.
\textit{Adv. Appl. Math.} \textbf{6} 4--22.

\bibitem{SJ12}
\textsc{Shivaswamy, P.} and \textsc{Joachims, T.} (2012).
Multi-armed bandit problems with history.
\textit{Proceedings of the Fifteenth International Conference on Artificial Intelligence and Statistics} 
\textbf{22} 1046--1054. 

\bibitem{SB98}
\textsc{Sutton, B.} and \textsc{Barto,, A.} (1998).
\textit{Reinforcement Learning, an Introduction}.
MIT Press, Cambridge.

\bibitem{TL10}
\textsc{Tekin, C.} and \textsc{Liu, M.} (2010>
Online algorithms for the multi-armed bandit problem with Markovian rewards.
\textit{48th Annual Allerton Conference on Communication, Control and Computing},
1675--1682.

\bibitem{TS85}
\textsc{Thathacher, V.} and \textsc{Sastry, P.S.} (1985).
A class of rapidly converging algorithms for learning automata.
\textit{IEEE Trans. Systems, Man Cyber.} \textbf{16} 168--175.

\bibitem{Tho33}
\textsc{Thompson, W.} (1933).
On the likelihood that one unknown probability exceeds another in view of the evidence of two samples.
\textit{Biometrika} \textbf{25} 285--294. 

\bibitem{YL91}
\textsc{Yakowitz, S.} and \textsc{Lowe, W.} (1991).
Nonparametric bandit problems.
\textit{Ann. Oper. Res.} \textbf{28} 297--312.
\end{thebibliography}
\end{document}